\documentclass[final,leqno,onefignum,onetabnum,sort&compress]{siamltex}
\usepackage{fancyhdr,adjustbox,circuitikz,pgfplots,comment,makecell,arydshln,csl}

\newcounter{book}
\usepackage{epsfig, graphicx}
\usepackage{makeidx}
\usepackage{multicol,multirow}
\usepackage{blkarray}
\usepackage{upgreek}
\usepackage{amsmath, amsfonts, amssymb, bbm, bm, mathtools, dsfont}
\usepackage{ulem}
\usepackage{crop}
\usepackage{xcolor,color}
\usepackage{todonotes}
\usepackage{hyperref,url}
\usepackage{subfig}
\usepackage[overload]{empheq}
\usepackage{cleveref,scalerel}
\usepackage{xspace,ifthen,mleftright,scalerel}
\usepackage{algorithmic}
\usepackage{xparse}
\usepackage{cite}

\usepackage{tikz}
\usepackage{tikz-qtree}
\usepackage{tikz-cd}

\usepackage{subfiles}
\usepackage{stmaryrd}
\usepackage{pdflscape}
\usepackage{breqn}
\usepackage{rotating}

\usepackage{array}
\newcolumntype{C}[1]{>{\centering\arraybackslash}m{#1}}  
\newcolumntype{L}[1]{>{\raggedleft\arraybackslash}m{#1}} 
\newcolumntype{R}[1]{>{\raggedleft\arraybackslash}m{#1}} 

\DeclareMathSizes{10}{10}{7}{4}
\DeclareMathSizes{11}{11}{8}{5}
\newcommand{\slow}{{\scalebox{0.7}{$\textsc{s}$}}}
\newcommand{\fast}{{\scalebox{0.7}{$\textsc{f}$}}}

\renewcommand{\slow}{{\textsc{s}}}
\renewcommand{\fast}{{\textsc{f}}}

\NewDocumentCommand{\component}{m m m}{^{ \left\{ \mkern-1mu #1 \mkern-1mu \right\} #2 \IfBooleanT{#3}{T} }}
\RenewDocumentCommand{\component}{m m m}{^{\left\{ \mkern-1mu \scalebox{0.65}{$#1$} \mkern-1mu \right\} #2 \IfBooleanT{#3}{T} }}
\NewDocumentCommand{\comp}{m O{} s}{\component{#1}{#2}{#3}}
\NewDocumentCommand{\F}{O{} s}{\component{\fast}{#1}{#2}}
\NewDocumentCommand{\FL}{O{\lambda} O{} s}{\component{\fast, \scalebox{0.8}{$#1$}}{#2}{#3}}
\RenewDocumentCommand{\S}{O{} s}{\component{\slow}{#1}{#2}}
\NewDocumentCommand{\SL}{O{\lambda} O{} s}{\component{\slow, \scalebox{0.8}{$#1$}}{#2}{#3}}
\NewDocumentCommand{\FF}{O{} s}{\component{\fast, \fast}{#1}{#2}}
\NewDocumentCommand{\FFL}{O{\lambda} O{} s}{\component{\fast, \fast, \scalebox{0.8}{$#1$}}{#2}{#3}}
\NewDocumentCommand{\FS}{O{} s}{\component{\fast, \slow}{#1}{#2}}
\NewDocumentCommand{\FSL}{O{\lambda} O{} s}{\component{\fast, \slow, \scalebox{0.8}{$#1$}}{#2}{#3}}
\NewDocumentCommand{\SF}{O{} s}{\component{\slow, \fast}{#1}{#2}}
\NewDocumentCommand{\SFL}{O{\lambda} O{} s}{\component{\slow, \fast, \scalebox{0.8}{$#1$}}{#2}{#3}}
\RenewDocumentCommand{\SS}{O{} s}{\component{\slow, \slow}{#1}{#2}}

\newcommand{\nvar}{d}
\newcommand{\nparts}{\mathrm{N}}

\newcommand{\fun}{\mathbf{f}}

\newcommand{\Lb}{\mathbf{L}}

\newcommand{\y}{\mathbf{y}}

\newcommand{\yf}[1][]{%
   \ifthenelse{ \equal{#1}{} }
      {{\y\F}}
      {{\y\F_{#1}}}
}
\newcommand{\dotyf}{{\dot{\y}\F}}

\newcommand{\ys}[1][]{%
   \ifthenelse{ \equal{#1}{} }
      {{\y\S}}
      {{\y\S_{#1}}}
}

\newcommand{\dotys}{{\dot{\y}\S}}

\newcommand{\mratio}{{\rm M}}

\newcommand{\lambdas}{{\lambda\S}}
\newcommand{\etas}{{\eta\S}}
\newcommand{\zs}{{z\S}}
\newcommand{\ws}{{w\S}}
\newcommand{\lambdaf}{{\lambda\F}}
\newcommand{\etaf}{{\eta\F}}
\newcommand{\zf}{{z\F}}
\newcommand{\wf}{{w\F}}

\newcommand{\Id}{\mathbf{I}}
\newcommand{\Idstage}{\Id_{s \times s}}
\newcommand{\Idvar}{{\Id_{\nvar \times \nvar}}} 
\newcommand{\one}{\boldsymbol{1}} 
\newcommand{\zero}{\boldsymbol{0}}
\newcommand{\sfrac}[2]{\mbox{\footnotesize$\displaystyle\frac{#1}{#2}$}} 

\makeatletter
\newcommand{\kron}[1]{{\,\mathrlap{\otimes}{\mathrlap{\hspace{0.275em}\tikz{\path[draw=white,fill=white] (0.025em,0.025em) rectangle (0.22em,0.38em);}}}\scalebox{0.4}{\raisebox{3.5pt}{\hspace{0.7em}#1}} \hspace{0.5em} }}
\makeatother

\def\Re{\mathds{R}}

\newtheorem{remark}{Remark}
\newtheorem{example}{Example}

\newcommand{\Vertiii}{\vert\kern-0.22ex\vert\kern-0.22ex\vert}
\newcommand{\Vertiv}{\vert\kern-0.22ex\vert\kern-0.22ex\vert\kern-0.22ex\vert}

\newcommand{\A}{{\mathbf{A}}}
\newcommand{\G}{{\mathbf{G}}}
\renewcommand{\b}{{\mathbf{b}}}
\newcommand{\B}{{\mathbf{B}}}
\renewcommand{\c}{{\mathbf{c}}}

\newcommand{\g}{{\mathbf{g}}}

\newenvironment{butchertableau}[2][1.3]
	{\def\arraystretch{#1}\array{#2}}
	{\endarray}

\newcommand{\Ys}[1][]{%
   \ifthenelse{ \equal{#1}{} }
      {{\mathbf{Y}_{\slow}}}
      {{\mathbf{Y}_{\slow,#1}}}
}
\newcommand{\Yf}[1][]{%
   \ifthenelse{ \equal{#1}{} }
      {{\mathbf{Y}_{\fast}}}
      {{\mathbf{Y}_{\fast,#1}}}
}

\renewcommand{\b}{\mathbf{b}\hspace{-1.4ex}\mathbf{b}}
\renewcommand{\c}{\mathbf{c}\hspace{-1.4ex}\mathbf{c}}
\renewcommand{\g}{\mathbf{g}\hspace{-1.4ex}\mathbf{g}}
\newcommand{\e}{\mathbf{e}\hspace{-1.4ex}\mathbf{e}}

\providecommand{\rosb}{\mathtt{b}}
\providecommand{\rosc}{\mathtt{c}}
\providecommand{\rosg}{\mathtt{g}}
\providecommand{\rose}{\mathtt{e}}
\providecommand{\rosalpha}{{\boldsymbol{\alpha}}}
\providecommand{\rosgamma}{{\boldsymbol{\gamma}}}
\providecommand{\rosbeta}{{\boldsymbol{\beta}}}

\allowdisplaybreaks

\begin{document}

\newcommand{\mytitle}{Multirate Linearly-Implicit GARK Schemes}

\ifnum \value{book}=1
	\csltitle{\mytitle}
	\cslauthor{Michael G\"unther and Adrian Sandu}
	\cslemail{guenther@uni-wuppertal.de, sandu@cs.vt.edu}
	\cslreportnumber{2}
	\cslyear{21}
	\csltitlepage
\fi

\title{\MakeUppercase{\mytitle}\thanks{
The work of Sandu was supported by awards NSF CCF--1613905. NSF ACI--1709727, NSF CDS\&E-MSS--1953113, and by the Computational Science Laboratory at Virginia Tech. The work of G\"unther was supported  by the European Union's Horizon 2020 research and innovation programme under the Marie Sklodowska-Curie Grant Agreement No.~765374, ROMSOC.}
}

\author{
Michael G\"unther\thanks{Bergische Universit\"at Wuppertal,
        Institute of Mathematical Modelling, Analysis and Computational
        Mathematics (IMACM), Gau\ss strasse 20, D-42119 Wuppertal, Germany 
        ({\tt guenther@uni-wuppertal.de})}
        \and
Adrian Sandu\thanks{Virginia Polytechnic Institute and State University, Computational Science Laboratory, Department of Computer Science, 2202 Kraft Drive, Blacksburg, VA 24060, USA ({\tt sandu@cs.vt.edu})}.
}

%

\maketitle

\begin{abstract}
Many complex applications require the solution of initial-value problems where some components change fast, while others vary slowly.  Multirate schemes apply different step sizes to resolve different components of the system, according to their dynamics, in order to achieve increased computational efficiency. The stiff components of the system, fast or slow, are best discretized with implicit base methods in order to ensure numerical stability. To this end, linearly implicit methods are particularly attractive as they solve only linear systems of equations at each step. 

This paper develops the Multirate GARK-ROS/ROW (MR-GARK-ROS/ROW) framework for linearly-implicit multirate time integration. The order conditions theory considers both exact and approximative Jacobians. The effectiveness of implicit multirate methods depends on the coupling between the slow and fast computations; an array of efficient coupling strategies and the resulting numerical schemes are analyzed. Multirate infinitesimal step linearly-implicit methods, that allow arbitrarily small micro-steps and offer extreme computational flexibility, are constructed. The new unifying framework includes existing multirate Rosenbrock(-W) methods as particular cases, and opens the possibility to develop new classes of highly effective linearly implicit multirate integrators. 

\end{abstract}



\pagestyle{myheadings}
\thispagestyle{plain}
\markboth{MICHAEL G\"UNTHER AND ADRIAN SANDU}{MULTIRATE GARK-ROS/ROW SCHEMES}


\section{Introduction}


Multiphysics applications lead to initial-value problems in  {\it additively} partitioned form  
\begin{equation} 
\label{eqn:additive-ode}
\y'= \fun(\y) = \sum_{m=1}^\nparts \fun^{\{m\}} (\y)\,, \quad  \y(t_0)=\y_0 \in \Re^{\nvar},
\end{equation}
where the right-hand side $\fun: \Re^{\nvar} \rightarrow \Re^{\nvar}$ is split into $\nparts$ different components 
based on, for example, stiffness (stiff/non-stiff), nonlinearity (linear/non-linear), dynamical behavior (fast/slow), and evaluation cost (cheap/expensive).  Additive partitioning also includes the special case of component partitioning where
the solution vector is split into disjoint sets~\cite{Sandu_2016_GARK-MR}.

Multimethods are effective numerical solvers for multiphysics and multiscale applications \eqref{eqn:additive-ode} that treat each component with an appropriate time discretization and time step, which are carefully coordinated such that the overall solution has the desired accuracy and stability properties. Examples of widely used multimethods include implicit-explicit \cite{Sandu_2015_IMEX-TSRK,Sandu_2016_highOrderIMEX-GLM,Sandu_2014_IMEX_GLM_Extrap ,Sandu_2014_IMEX-GLM,Sandu_2014_IMEX-RK,Kennedy_2019_IMEX-ARK} and multirate  schemes  \cite{Gear_1984_MR-LMM,Rice_1960_splitRK,Bremicker_2017_MGARK}. 

This paper focuses on partitioned systems \eqref{eqn:additive-ode} where some components are fast varying and have small evaluation costs, while other are slowly varying but their evaluation is expensive. Multirate schemes exploit this different dynamical behavior by applying small step sizes to the fast part and large step sizes to the slow part, while ensuring the order and stability of the overall numerical integration scheme. Significant computational savings are made possible by the less frequent evaluation of the slow expensive components. Multirate schemes have been developed using base methods such as Runge-Kutta \cite{Sandu_2007_MR_RK2,Rice_1960_splitRK,Andrus_1979_MR,Andrus_1993_MR-stability,Kvaerno_2000_stability-MRK,Kvaerno_1999_MR-RK}, Linear Multistep \cite{Sandu_2009_MR_LMM,Gear_1984_MR-LMM,Kato_1999}, Rosenbrock~\cite{Guenther_2001_MR-PRK,Guenther_1993_MR-ROW,Savcenco_2007_MR-Rosenbrock}, and extrapolation \cite{Engstler_1997_MR-extrapolation,Sandu_2008_MR-EXTRAP,Sandu_2013_extrapolatedMR,Sandu_2009_ICNAAM-Multirate}.

The General-structure Additive Runge-Kutta (GARK) formalism introduced in \cite{Sandu_2015_GARK} defines a comprehensive framework for studying a large class of partitioned Runge-Kutta based schemes  for solving~\eqref{eqn:additive-ode}. It allows for different stage values with different components of the right hand side. Though in principle equivalent to Additive Runge-Kutta schemes~\cite{Kennedy_2003_ARK},``the advantage of the GARK formulation is that it clarifies the coupling between the various methods, in addition to eliminating zero quadrature weights in the ARK formalism, hence the analysis of special cases''~\cite{zanna2020}. Examples of partitioned methods developed in the GARK framework include~\cite{Sandu_2020_GARK-adjoint,Sandu_2016_GARK-MR,zanna2020}.
MR-GARK and MRI-GARK frameworks~\cite{Sandu_2016_GARK-MR,Sandu_2019_MR-GARK_High-Order,Sandu_2019_MRI-GARK,Sandu_2020_MRI-GARK_Coupled,Sandu_2021_MR-DAE,Sandu_2021_MR-GARK_Implicit} define general classes of multirate Runge-Kutta methods based on the GARK formalism. 

The stiff components of the system \eqref{eqn:additive-ode} are best discretized with implicit base methods in order to ensure numerical stability. Linearly implicit methods enjoy the same stability properties as the implicit schemes, but solve only linear systems of equations at each step. Rosenbrock-Wanner (ROS) methods \cite{Rosenbrock_1963,Hairer_book_II} are linearly implicit Runge--Kutta schemes that use the exact Jacobian of the right hand side function in the computational process;  Rosenbrock-W methods \cite{Steihaug_1979} allow arbitrary approximations of the Jacobian.
In a recent paper~\cite{Sandu_2021_GARK-ROS} the authors have generalized the GARK approach to partitioned linearly-implicit schemes based on using exact (GARK-ROS) and inexact (GARK-ROW) Jacobian information. 

This paper proposes the Multirate GARK-ROS/ROW (MR-GARK-ROS/ROW) framework for linearly-implicit multirate time integration. The new elements developed herein are as follows. 
MR-GARK-ROS/ROW provides a unifying formalism for linearly-implicit multirate schemes of Runge-Kutta type, that includes as particular cases existing methods such as ~\cite{Guenther_2001_MR-PRK,Guenther_1993_MR-ROW,Savcenco_2007_MR-Rosenbrock}.
A general order conditions theory is developed for methods with exact and with inexact Jacobian information. 
The effectiveness of implicit multirate methods hinges upon the manner in which the slow and fast system computations are coupled; to this end, an array of various efficient coupling approaches is proposed.
Two particular coupling structures, which hold considerable promise for practical applications, are studied in detail: compound-first-step, where the macro-step is coupled with the first micro-step only, and step-predictor-corrector, where a macro-step carried out over the entire system is followed by a recalculation of the fast components using small steps.
Finally, multirate infinitesimal step (MRI) methods \cite{Wensch_2009_MIS,Knoth_1998_MR-IMEX,Sandu_2019_MRI-GARK,Sandu_2020_MRI-GARK_Coupled} that allow arbitrarily small micro-steps offer extreme computational flexibility, since the fast subsystem can be solved with any discretization method and sequence of step sizes. The work develops general MRI-GARK ROS/ROW schemes,  as well as MRI step-predictor-corrector methods; order condition theories are provided for both families.

The paper is organized as follows. \Cref{sec:review-GARK} reviews basic aspects of GARK-ROS/ROW methods. 
The new MR-GARK ROS/ROW formalism is defined in \Cref{sec:define-MGARK-ROS/ROW}. Using the partitioned GARK ROS/ROW framework, general order conditions are derived in \Cref{sec:order-conditions}.
A linear stability analysis is performed in \Cref{sec:linear-stability}.  
Various slow-fast coupling strategies for an efficient computation are discussed in~\Cref{sec:coupling-strategies}. 
Compound-first-step schemes are analyzed in \Cref{sec:compound-first-step}, and step-predictor-corrector (SPC) methods  in \Cref{sec:Step-predictor-corrector}. 
Multirate infinitesimal step  (MRI-GARK ROS/ROW) schemes are developed in \Cref{sec:infinitesimal-step}, and multirate infinitesimal step SPC methods  in \Cref{sec:MRI-SPC-methods}.
Conclusions and an outlook of future work are given in \Cref{sec:conclusions}.

\section{Linearly implicit GARK schemes}
\label{sec:review-GARK}
\ifnum \value{book}=1
In the following we recapitulate some basic facts on GARK \cite{Sandu_2015_GARK},  Multirate GARK GARK-ROS/ROW \cite{Sandu_2021_GARK-ROS} schemes.
\fi 

\ifnum \value{book}=1
%
\subsection{GARK  schemes}

\begin{definition}[GARK methods~\cite{Sandu_2015_GARK}] 
\label{def:GARK-different-stages}
One step of a GARK scheme with a $\nparts$-way partitioning of the right hand side  \eqref{eqn:additive-ode} reads:
\begin{subequations}
\label{eqn:GARK}
\begin{eqnarray}
\label{eqn:GARK-stage}
\mathbf{k}^{\{q\}} &=& h\,\fun^{\{q\}}\left( 
\one_s \otimes \y_{n-1}  + \sum_{m=1}^\nparts \A^{\{q,m\}} \kron{\nvar}  \mathbf{k}^{\{m\}} \right),    
\quad q = 1,\dots,\nparts, \\ 
\label{eqn:GARK-solution}
\y_{n} &=& \y_{n-1} + \sum_{m=1}^\nparts  \b^{\{q\}}\,\!^T  \kron{\nvar}  \mathbf{k}^{\{q\}}.
\end{eqnarray}
\end{subequations} 
Here $\one_s \in \Re^s$ is a vector of ones, $\otimes$ is the Kronecker product, and the following matrix notation is used:
\begin{subequations}
\label{eqn:GARK-matrix-notation}
\begin{equation}
\label{eqn:K-matrix-notation}
\mathbf{k}^{\{q\}} \coloneqq 
\begin{bmatrix}
\mathbf{k}_1^{\{q\}} \\ 
\vdots \\ \mathbf{k}_{s^{\{q\}}}^{\{q\}}
\end{bmatrix} \in \Re^{\nvar s^{\{q\}}},
\qquad
\begin{array}{c}
\mathbb{\upalpha}^{\{m,q\}} \kron{\nvar}  \mathbf{k}^{\{q\}} \coloneqq \left( \mathbb{\upalpha}^{\{m,q\}} \otimes \Idvar \right)\, \mathbf{k}^{\{q\}}, \\
\mathbb{\upgamma}^{\{m,q\}} \kron{\nvar}  \mathbf{k}^{\{q\}} \coloneqq \left( \mathbb{\upgamma}^{\{m,q\}} \otimes \Idvar \right)\, \mathbf{k}^{\{q\}},
\end{array}
\end{equation}
with $\A^{\{m,q\}} \in \Re^{s^{\{m\}} \times s^{\{q\}}}$, $\Idvar \in \Re^{\nvar \times \nvar}$ the identity matrix, and
\begin{equation}
\label{eqn:Fun-matrix-notation}
\begin{split}
& \fun^{\{q\}}\left( 
\one_s \otimes \y_{n-1}  + \sum_{m=1}^\nparts \A^{\{q,m\}} \kron{\nvar}  \mathbf{k}^{\{m\}} \right) \coloneqq \\
& \qquad \begin{bmatrix} \displaystyle \fun^{\{q\}}( \y_{n-1}  + \sum_{m=1}^\nparts \sum_{j=1}^{s^{\{m\}}} a_{i,j}^{\{q,m\}}\, \mathbf{k}_j^{\{m\}})\end{bmatrix}_{i \le i \le s^{\{q\}}}.
\end{split}
\end{equation}
\end{subequations}
\end{definition}
The corresponding generalized Butcher tableau is
\begin{equation}
\label{eqn:general-Butcher-tableau-different-stages}
\renewcommand{\arraystretch}{1.25} 
\begin{array}{ccc}
\A^{\{1,1\}} &  \ldots & \A^{\{1,\nparts\}} \\
\vdots &  & \vdots \\
\A^{\{\nparts,1\}} &  \ldots & \A^{\{\nparts,\nparts\}} \\ \Xhline{1.5pt}
\b^{\{1\}} &  \ldots &\b^{\{\nparts\}}
\end{array}.
%
\end{equation}

\begin{remark}
In contrast to traditional additive methods~\cite{Kennedy_2003_ARK,Kennedy_2019_IMEX-ARK} different stage values are used with different components of the right hand side. 
The methods $(\A^{\{q,q\}},\b^{\{q\}})$ in \eqref{eqn:GARK} can be regarded as stand-alone integration schemes
applied to each individual component $q$. The off-diagonal matrices $\A^{\{q,m\}}$, $m \ne q$, can be 
viewed as a coupling mechanism among components. 
\end{remark}

\begin{definition}[Internally consistent GARK methods]
A GARK scheme \eqref{eqn:GARK} is called internally consistent if
\begin{equation}
\label{eqn:GARK-internal-consistency}
\sum_{j=1}^{s^{\{1\}}} a_{i,j}^{\{q,1\}} 
= \dots = \sum_{j=1}^{s^{\{\nparts\}}} a_{i,j}^{\{q,\nparts\}} = c_i^{\{q\}}\,,
\quad i=1,\dots, s^{\{q\}}\,, \quad q=1,\dots,\mratio.
\end{equation}
\end{definition}
The internal consistency condition \eqref{eqn:GARK-internal-consistency} ensures that all components of the 
stage vectors are calculated at the same internal approximation times.

\subsection{Multirate GARK  schemes}
Multirate GARK schemes \cite{Sandu_2016_GARK-MR} seek to exploit the multiscale behavior given in different dynamics.  Consider a two-way additively partitioned system \eqref{eqn:additive-ode} of the form
\begin{equation} 
\label{eqn:multirate-ode}
\y'= \fun(\y) = \fun\S(\y) + \fun\F(\y)\,, \quad  \y(t_0)=\y_0 \in \Re^{\nvar},
\end{equation}
with a slow part $\fun\S$ that is computationally expensive and a fast part $\fun\F$ that is inexpensive to evaluate. 
Solving  the slow  expensive component with a large step $H$, and the fast inexpensive one with small steps $h=H/\mratio$, allows for reducing the overall computational cost. A multirate generalization of \eqref{eqn:GARK} with $\mratio$ micro steps $h=H/\mratio$  is then given by the following.

\begin{definition}[Multirate GARK method~\cite{Sandu_2016_GARK-MR}] One macro-step of a generalized additive multirate Runge-Kutta method with $\mratio$ equal micro-steps reads:
\begin{subequations}
\label{eqn:MGARK}
\begin{eqnarray}
\mathbf{k}\FL[\lambda] &=& 
h\,\fun\F\left( 
\one_s \otimes \y_{n-1}  + 
\sum_{\ell=1}^{\lambda-1} \b\F* \kron{\nvar} \mathbf{k}\FL[\ell] +
\A\FF \kron{\nvar} \mathbf{k}\FL[\lambda]   \right. \nonumber \\
& & \qquad\qquad +  \A\FSL[\lambda] \kron{\nvar} \mathbf{k}\S
\Bigg), \quad {\rm for}\quad \lambda=1,\ldots,\mratio; \\
\mathbf{k}\S &=& H\,\fun\S\left( 
\one_s \otimes \y_{n-1}  + \sum_{\lambda=1}^{\mratio}
 \A\SFL[\lambda] \kron{\nvar}  \mathbf{k}\FL[\lambda] 
+ \A\SS \kron{\nvar}  \mathbf{k}\S
\right); \\
\y_{n} &=& \y_{n-1} + \sum_{\lambda=1}^{\mratio} \b\F*\kron{\nvar}  \mathbf{k}\FL[\lambda] + \b\S*\kron{\nvar}  \mathbf{k}\S.
\end{eqnarray}
\end{subequations} 

The base  schemes are Runge-Kutta methods, {$(A\SS,b\S)$} for the slow component and {$(A\FF,b\F)$} for the fast component. The coefficients $A\SFL[\lambda]$, $A\FSL[\lambda]$ realize the coupling between the two components.
\end{definition}
The method \eqref{eqn:GARK} can be
written as a GARK scheme~\eqref{eqn:GARK}
over the macro-step $H$ . 
The corresponding Butcher tableau 
\eqref{eqn:general-Butcher-tableau-different-stages}
reads
\begin{equation}
\label{eqn:mrRK-butcher}
\renewcommand{\arraystretch}{1.5}
\begin{array}{c|c}
{\A}\FF & {\A}\FS  \\ \hline
{\A}\SF &{\A}\SS \\ \Xhline{1.5pt}
{\b}\F  \,\!^T & \b\S  \,\!^T
\end{array} 
~~  \raisebox{-9pt}{$\coloneqq$} ~~ 
\raisebox{28.5pt}{$
\begin{array}{cccc|cccc}  
\sfrac{1}{\mratio} \A\FF      &          0                   & \cdots & 0 & \A\FSL[1]  \\
\sfrac{1}{\mratio} \one  \b\F* & \sfrac{1}{\mratio} \A\FF        & \cdots & 0 &  \A\FSL[2]  \\
\vdots                     &                             & \ddots &   & \vdots  \\
\sfrac{1}{\mratio} \one  \b\F* & \sfrac{1}{\mratio} \one  \b\F*   & \ldots & \sfrac{1}{\mratio} \A\FF & \A\FSL[\mratio] \\
\hline 
\sfrac{1}{\mratio} \A\SFL[1] & \sfrac{1}{\mratio} \A\SFL[2] & \cdots & \sfrac{1}{\mratio} \A\SFL[\mratio] & \A\SS   \\   \Xhline{1.5pt}
\sfrac{1}{\mratio} \b\F* & \sfrac{1}{\mratio} \b\F* & \ldots & \sfrac{1}{\mratio} \b\F* & \b\S*  
\end{array}
$}.
\renewcommand{\arraystretch}{1.0}
\end{equation}
\fi

\subsection{GARK-ROS and GARK-ROW schemes}

The class of linearly-implicit GARK-ROS and GARK-ROW schemes was developed in ~\cite{Sandu_2021_GARK-ROS} , in analogy to the extension of (explicit) Runge-Kutta schemes to Rosenbrock-Wanner (ROS) and Rosenbrock-W  (ROW) schemes. One GARK-ROS/ROW step applied to \eqref{eqn:additive-ode} reads:
\begin{subequations}
\label{eqn:GARK-ROS/ROW}
\begin{eqnarray}
\label{eqn:GARK-ROS/ROW-stage}
\mathbf{k}^{\{q\}} &=& h\,\fun^{\{q\}}\left( 
\one_s \otimes \y_{n-1}  + \sum_{m=1}^\nparts \pmb{\upalpha}^{\{q,m\}} \kron{\nvar}  \mathbf{k}^{\{m\}} \right)    \\
\nonumber
&& + \bigl( \Idstage \otimes h\,\Lb^{\{q\}} \bigr)\,\sum_{m=1}^\nparts  \pmb{\upgamma}^{\{q,m\}} \kron{\nvar}  \mathbf{k}^{\{m\}}, \quad q = 1,\dots,\nparts, \\ 
\label{GARK-ROS/ROW-solution}
\y_{n} &=& \y_{n-1} + \sum_{m=1}^\nparts  \b^{\{m\}}\,\!^T  \kron{\nvar}  \mathbf{k}^{\{m\}}.
\end{eqnarray}
\end{subequations} 
\ifnum \value{book}=1
where we use the matrix notation \eqref{eqn:GARK-matrix-notation}. 
\else
Here $\one_s \in \Re^s$ is a vector of ones, $\Idstage \in \Re^{s \times s}$ is the identity matrix, $\otimes$ is the Kronecker product, and the following matrix notation is used:
\begin{subequations}
\label{eqn:GARK-matrix-notation}
\begin{equation}
\label{eqn:K-matrix-notation}
\mathbf{k}^{\{q\}} \coloneqq 
\begin{bmatrix}
{\scriptstyle \mathbf{k}_1^{\{q\}} } \\ 
{\scriptstyle \vdots} \\ {\scriptstyle \mathbf{k}_{s^{\{q\}}}^{\{q\}} }
\end{bmatrix} \in \Re^{\nvar s^{\{q\}}},
\qquad
\begin{array}{c}
\pmb{\upalpha}^{\{m,q\}} \kron{\nvar}  \mathbf{k}^{\{q\}} \coloneqq \left( \pmb{\upalpha}^{\{m,q\}} \otimes \Idvar \right)\, \mathbf{k}^{\{q\}}, \\[3pt]
\pmb{\upgamma}^{\{m,q\}} \kron{\nvar}  \mathbf{k}^{\{q\}} \coloneqq \left( \pmb{\upgamma}^{\{m,q\}} \otimes \Idvar \right)\, \mathbf{k}^{\{q\}},
\end{array}
\end{equation}
with $\pmb{\upalpha}^{\{m,q\}},  \pmb{\upgamma}^{\{m,q\}}\in \Re^{s^{\{m\}} \times s^{\{q\}}}$ and
\begin{equation}
\label{eqn:Fun-matrix-notation}
\fun^{\{q\}}\left( 
{\scriptstyle
\one_s \otimes \y_{n-1}  + \sum_{m=1}^\nparts \pmb{\upalpha}^{\{q,m\}} \kron{\nvar}  \mathbf{k}^{\{m\}} }\right) 
\coloneqq 
\begin{bmatrix}
\fun^{\{q\}}\bigl( {\scriptstyle \y_{n-1}  + \sum_{m=1}^\nparts \sum_{j=1}^{s^{\{m\}}} \upalpha_{1,j}^{\{q,m\}}\, \mathbf{k}_j^{\{m\}} }\bigr) \\  
\ifnum \value{book}=1
\fun^{\{q\}}\bigl( \y_{n-1}  + \sum_{m=1}^\nparts \sum_{j=1}^{s^{\{m\}}} \upalpha_{2,j}^{\{q,m\}}\, \mathbf{k}_j^{\{m\}}\bigr)\\ 
\fi
{\scriptstyle \vdots} \\ 
\fun^{\{q\}}\bigl({\scriptstyle \y_{n-1}  + \sum_{m=1}^\nparts \sum_{j=1}^{s^{\{m\}}} \upalpha_{s^{\{q\}},j}^{\{q,m\}}\, \mathbf{k}_j^{\{m\}} }\bigr)
\end{bmatrix}.
\end{equation}
\end{subequations}
\fi

The matrices $\pmb{\upalpha}^{\{q,m\}}$ are strictly lower triangular  and $\pmb{\upgamma}^{\{q,m\}}$ are lower 
triangular. Depending on the choice of matrices $\Lb^{\{q\}} \in \Re^{\nvar \times \nvar}$ one distinguishes several types of methods, as follows:
\begin{itemize}
    \item  GARK-ROS schemes use the exact Jacobian information, i.e., $\Lb^{\{q\}} \coloneqq \fun_{\y}^{\{q\}}(\y_{n-1})$ are the Jacobians of the component functions evaluated at the current solution;
    \item GARK-ROW schemes allow any approximation of the Jacobian, i.e., $\Lb^{\{q\}}$ may be arbitrary; 
    \item In the case of GARK-ROS schemes with time-lagged Jacobians one has $\Lb^{\{q\}}=\fun_{\y}^{\{q\}}(\y_{n-1}) + \mathcal{O}(h)$.
\end{itemize}

The scheme \eqref{eqn:GARK-ROS/ROW} is characterized by the extended Butcher tableau (with $\pmb{\upbeta}^{\{q,m\}} \coloneqq \pmb{\upalpha}^{\{q,m\}}+\pmb{\upgamma}^{\{q,m\}}$)
\begin{equation}
\label{eqn:GARK-Rosenbrock-butcher}
\begin{butchertableau}{c|c}
\A  & \G  \\ 
\Xhline{1.5pt}
\b^T &
\end{butchertableau}
~=~
\raisebox{17pt}{$
\begin{butchertableau}{ccc | ccc}
\pmb{\upalpha}^{\{1,1\}} &  \cdots & \pmb{\upalpha}^{\{1,\nparts\}} & \pmb{\upgamma}^{\{1,1\}} &  \ldots & \pmb{\upgamma}^{\{1,\nparts\}}  \\
\vdots &  \ddots & \vdots & \vdots & \ddots & \vdots \\
\pmb{\upalpha}^{\{\nparts,1\}} & \ \cdots & \pmb{\upalpha}^{\{\nparts,\nparts\}} & \pmb{\upgamma}^{\{\nparts,1\}} &  \ldots & \pmb{\upgamma}^{\{\nparts,\nparts\}} \\ 
\Xhline{1.5pt}
\b^{\{1\}T} &  \cdots &\b^{\{\nparts\}T} &  &    &
\end{butchertableau}
$}.
\end{equation}

\begin{remark}[GARK-ROS and GARK-ROW scheme structure]
\begin{itemize}
\item Similar to GARK, the scheme \eqref{eqn:GARK-ROS/ROW} uses only one function $\fun^{\{q\}}$ evaluation for the increment $\mathbf{k}^{\{q\}}$, and linear combinations of increments $\mathbf{k}^{\{m\}}$  both as function arguments and as additive terms. 
\item For all $\pmb{\upgamma}^{\{q,m\}}=0$ \eqref{eqn:GARK-ROS/ROW} is an explicit GARK scheme.
\item If $\pmb{\upgamma}^{\{q,m\}}=0$ for all $m>q$ then increments can be computed sequentially in the following order: $k_1^{\{1\}}, \ldots, k_1^{\{\nparts\}}, k_2^{\{1\}},\ldots, k_{s^{\{\nparts\}}}^{\{\nparts\}}$.
\end{itemize}
\end{remark}

\begin{theorem}[GARK-ROS order conditions~\cite{Sandu_2021_GARK-ROS}]
The GARK-ROS order conditions \cref{eqn:GARK-ROS/ROW}  are the same as the Rosenbrock order conditions  \cite{Hairer_book_II}, except that the method coefficients are also labelled according to partition indices. In the order conditions, in each sequence of matrix multiplies, the partition indices are compatible according to matrix multiplication rules.
\end{theorem}

Let $\one^{\{n\}} \in \Re^{s^{\{n\}}}$ be a vector of ones. For brevity of notation we define:
\begin{equation}
\label{eqn:terminal-vectors}
\begin{split}
\c^{\{m,n\}} &\coloneqq \pmb{\upalpha}^{\{m,n\}} \,  \one^{\{n\}}, \quad
\g^{\{m,n\}} \coloneqq \pmb{\upgamma}^{\{m,n\}} \,  \one^{\{n\}}, \\
\e^{\{m,n\}} &\coloneqq \pmb{\upbeta}^{\{m,n\}} \, \one^{\{n\}}
= \c^{\{m,n\}} + \g^{\{m,n\}}.
\end{split}
\end{equation}
The GARK-ROS order four conditions read: 
\begin{subequations}
\label{eqn:GARK-ROS-order-conditions}
\begin{align}
\label{eqn:GARK-ROS-order1-conditions}
\textnormal{order 1:   } &\begin{cases}
\b^{\{m\}\,T} \, \one^{\{m\}} = 1, \end{cases} \textnormal{for } m = 1,\dots,\nparts; \\
\label{eqn:GARK-ROS-order2-conditions}
\textnormal{order 2:   } &
\begin{cases}
\b^{\{m\}\,T} \, \e^{\{m,n\}}  = \frac{1}{2}, \end{cases}  \textnormal{for } m,n = 1,\dots,\nparts; \\
\label{eqn:GARK-ROS-order3-conditions}
\textnormal{order 3:   } &
\begin{cases}
\b^{\{m\}\,T} \, \bigl( \c^{\{m,n\}} \times \c^{\{m,p\}} \bigr) = \frac{1}{3}, \\
\b^{\{m\}\,T} \, \pmb{\upbeta}^{\{m,n\}}\, \e^{\{n,p\}} = \frac{1}{6},
\end{cases}  \textnormal{for } m,n,p = 1,\dots,\nparts; \\
\label{eqn:GARK-ROS-order4-conditions}
\textnormal{order 4:   } &
\begin{cases}
\b^{\{m\}\,T} \,\big( \c^{\{m,n\}} \times  \c^{\{m,p\}}  \times  \c^{\{m,q\}} \big) = \frac{1}{4}, \\
\b^{\{m\}\,T} \, \Big((\pmb{\upalpha}^{\{m,n\}}\,\e^{\{n,p\}}) \times \c^{\{m,q\}} \Big) = \frac{1}{8}, \\
\b^{\{m\}\,T} \, \pmb{\upbeta}^{\{m,n\}}\, \left( \c^{\{n,p\}} \times \c^{\{n,q\}}  \right) = \sfrac{1}{12}, \\ 
\b^{\{m\}\,T} \,\pmb{\upbeta}^{\{m,n\}}\,\pmb{\upbeta}^{\{n,p\}}\,\e^{\{p,q\}} =  \sfrac{1}{24}, \\
 \textnormal{for } m,n,p,q = 1,\dots,\nparts.
\end{cases} 
\end{align}
\end{subequations}
%
%
\begin{theorem}[GARK-ROW order conditions~\cite{Sandu_2021_GARK-ROS}]
The GARK-ROW order conditions \cref{eqn:GARK-ROS/ROW}  are the same as the Rosenbrock-W order conditions \cite{Hairer_book_II}, except that the method coefficients are also labelled according to partition indices. In the order conditions, in each sequence of matrix multiplies, the partition indices are compatible according to matrix multiplication rules.
\end{theorem}

The GARK-ROW order four conditions read:
\begin{subequations}
\label{eqn:GARK-ROW-order-conditions}
\begin{alignat}{3}
\label{eqn:GARK-ROW-order1-conditions}
\quad\textnormal{   order 1:   } &
\begin{cases} \b^{\{m\}\,T} \, \one^{\{m\}} = 1, & \textnormal{for  } m = 1,\dots,\nparts; \end{cases} \\
\label{eqn:GARK-ROW-order2-conditions}
\textnormal{   order 2:   } &
\begin{cases}
\b^{\{m\}\,T} \, \c^{\{m,n\}}  = \sfrac{1}{2}, & \\
\b^{\{m\}\,T} \, \g^{\{m,n\}}  = 0,
& \textnormal{for  } m,n = 1,\dots,\nparts;
\end{cases}  \\
\label{eqn:GARK-ROW-order3-conditions}
\textnormal{   order 3:   } &
\begin{cases}
\b^{\{m\}\,T} \, \bigl( \c^{\{m,n\}} \times \c^{\{m,p\}} \bigr) = \sfrac{1}{3}, &
\b^{\{m\}\,T} \, \pmb{\upalpha}^{\{m,n\}}\, \c^{\{n,p\}} = \sfrac{1}{6}, \\
\b^{\{m\}\,T} \, \pmb{\upgamma}^{\{m,n\}}\,\c^{\{n,p\}} = 0, &
\b^{\{m\}\,T} \, \pmb{\upalpha}^{\{m,n\}}\,\g^{\{n,p\}}= 0, \\
\b^{\{m\}\,T} \, \pmb{\upgamma}^{\{m,n\}}\, \g^{\{n,p\}} = 0, &
\textnormal{for  } m,n,p = 1,\dots,\nparts;
\end{cases}
\end{alignat}
\begin{equation}
\label{eqn:GARK-ROW-order4-conditions}
\begin{split}
&\textnormal{order 4:   } \\
&\begin{cases}
\b^{\{m\}\,T} \, \left( \c^{\{m,n\}} \times \c^{\{m,p\}} \times \c^{\{m,q\}}\right) = \sfrac{1}{4}, \\
\b^{\{m\}\,T} \, \big((\pmb{\upalpha}^{\{m,n\}}\,\c^{\{n,p\}}) \times \c^{\{m,q\}} \big) = \sfrac{1}{8},  \\
\b^{\{m\}\,T} \, \pmb{\upalpha}^{\{m,n\}}\, (\c^{\{n,p\}} \times \c^{\{n,q\}}) = \sfrac{1}{12}, \\
\b^{\{m\}\,T}\,\pmb{\upalpha}^{\{m,n\}}\,\pmb{\upalpha}^{\{n,p\}}\, \c^{\{p,q\}} = \sfrac{1}{24}, \\ 
\b^{\{m\}\,T} \, \big((\pmb{\upalpha}^{\{m,n\}}\,\g^{\{n,p\}}) \times \c^{\{m,q\}} \big) = 0, \\
\b^{\{m\}\,T}\, \pmb{\upgamma}^{\{m,n\}}\, (\c^{\{n,p\}} \times \c^{\{n,q\}}) = 0, \\ 
\b^{\{m\}\,T} \,\pmb{\upgamma}^{\{m,n\}}\,\pmb{\upalpha}^{\{n,p\}}\, \c^{\{p,q\}} = 0, 
~~ \b^{\{m\}\,T}\,\pmb{\upalpha}^{\{m,n\}}\,\pmb{\upgamma}^{\{n,p\}}\,\c^{\{p,q\}} = 0, \\ 
\b^{\{m\}\,T}\,\pmb{\upalpha}^{\{m,n\}}\,\pmb{\upalpha}^{\{n,p\}}\,\g^{\{p,q\}} = 0, 
~~ \b^{\{m\}\,T}\,\pmb{\upgamma}^{\{m,n\}}\,\pmb{\upalpha}^{\{n,p\}}\,\g^{\{p,q\}} = 0, \\ 
\b^{\{m\}\,T}\,\pmb{\upalpha}^{\{m,n\}}\,\pmb{\upgamma}^{\{n,p\}}\, \g^{\{p,q\}} = 0,  
~~ \b^{\{m\}\,T}\,\pmb{\upgamma}^{\{m,n\}}\,\pmb{\upgamma}^{\{n,p\}}\,\c^{\{p,q\}} = 0, \\ 
\b^{\{m\}\,T}\,\pmb{\upgamma}^{\{m,n\}}\,\pmb{\upgamma}^{\{n,p\}}\, \g^{\{p,q\}} = 0, 
\quad \textnormal{for  } m,n,p,q = 1,\dots,\nparts.
\end{cases}
\end{split}
\end{equation}
\end{subequations}

\begin{remark}[Internal consistency~\cite{Sandu_2021_GARK-ROS}]
The order conditions~\eqref{eqn:GARK-ROS-order-conditions} and \eqref{eqn:GARK-ROW-order-conditions} simplify considerably for internally consistent GARK-ROS/ROW schemes, for which the vectors \eqref{eqn:terminal-vectors} satisfy:
\begin{equation}
\label{eqn:ROW-internal-consistency}
\c^{\{m,n\}} = \c^{\{m\}}, \qquad 
\g^{\{m,n\}} = \g^{\{m\}}, \qquad 
 \forall\, m,n = 1,\dots,\nparts.
\end{equation}
\end{remark}

\section{Multirate GARK ROS/ROW methods}
\label{sec:define-MGARK-ROS/ROW}
We aim at employing the GARK ROS/ROW formalism \eqref{eqn:GARK-ROS/ROW} to develop a class of multirate linearly-implicit schemes. 
\ifnum \value{book}=1
The philosophy is similar to the one used for developing Multirate GARK schemes \eqref{eqn:MGARK} 
based on the GARK framework \eqref{eqn:GARK} \cite{Sandu_2016_GARK-MR,Sandu_2015_GARK}. 
\fi

\subsection{Multirate GARK-ROS/ROW for additive splitting}
\ifnum \value{book}=0
Multirate GARK schemes \cite{Sandu_2016_GARK-MR} seek to exploit the multiscale behavior given in different dynamics.  Consider a two-way additively partitioned system \eqref{eqn:additive-ode} of the form
\begin{equation} 
\label{eqn:multirate-ode}
\y'= \fun(\y) = \fun\S(\y) + \fun\F(\y)\,, \quad  \y(t_0)=\y_0 \in \Re^{\nvar},
\end{equation}
with a slow part $\fun\S$ that is computationally expensive and a fast part $\fun\F$ that is inexpensive to evaluate. 
Solving  the slow  expensive component with a large macro-step $H$, and the fast inexpensive one with $\mratio$ micro-steps $h=H/\mratio$, allows for reducing the overall computational cost.  
\fi

\begin{definition}[Multirate GARK-ROS/ROW method] 
\label{def:MGARK-ROS/ROW}
One step of a multirate GARK-ROS/ROW method (MR-GARK-ROS/ROW for short) applied to \eqref{eqn:multirate-ode} computes the solution as follows.
The slow component is discretized with a ROS/ROW method $\bigl(  \rosb\S,\rosalpha\SS,\rosgamma\SS \bigr)$ and macro-step $H$. The fast component uses $\mratio$ micro-steps $h = H/\mratio$, and at each micro-step $\lambda$ a (possibly different) ROS/ROW method $\bigl(  \rosb\FL[\lambda],\rosalpha\FFL[\lambda],\rosgamma\FFL[\lambda] \bigr)$  is applied. 
The computational process reads:
%
%
\begin{subequations}
\label{eqn:MROS/MROW}
\begin{eqnarray}
\mathbf{k}\FL[\lambda] &=& h\,\fun\F\Bigl( \one_{s\F} \otimes \y_{n-1}  + 
\sum_{\ell=1}^{\lambda-1} \one\F\, \rosb\FL[\ell]*  \kron{\nvar} \mathbf{k}\FL[\ell] +
 \rosalpha\FFL[\lambda]  \kron{\nvar} \mathbf{k}\FL[\lambda]  \nonumber \\
 & & \nonumber \qquad + 
 \rosalpha\FSL[\lambda]  \kron{\nvar} \mathbf{k}\S \Bigr) \\
\label{eqn:MROS/MROW-KF}
&& + \Bigl( \Id_{s\F \times s\F} \otimes h\,\Lb\F \Bigr) \,
\Bigl( 
\rosgamma\FFL[\lambda]  \kron{\nvar} \mathbf{k}\FL[\lambda]
+ \rosgamma\FSL[\lambda] \kron{\nvar} \mathbf{k}\S \Bigr),  \\[3pt]
& &\quad  \textnormal{for}\quad \lambda=1,\ldots,\mratio;  \nonumber \\
\mathbf{k}\S &=& H\,\fun\S\left( \one_{s\S} \otimes \y_{n-1}  + 
\sum_{\lambda=1}^\mratio 
\rosalpha\SFL[\lambda] \kron{\nvar} \mathbf{k}\FL[\lambda] 
+  \rosalpha\SS  \kron{\nvar} \mathbf{k}\S \right)  \nonumber \\
\label{eqn:MROS/MROW-KS}
& &  +  \bigl( \Id_{s\S \times s\S} \otimes H\,\Lb\S) \,
\left(
\sum_{\lambda=1}^\mratio
\rosgamma\SFL[\lambda]  \kron{\nvar}  \mathbf{k}\FL[\lambda]
+ \rosgamma\SS  \kron{\nvar}  \mathbf{k}\S  \right); \\
\label{eqn:MROS/MROW-Sol}
\y_{n} &=& \y_{n-1} +
\sum_{\lambda=1}^{\mratio} \rosb\FL[\lambda]* \kron{\nvar} \mathbf{k}\FL[\lambda] +  \rosb\S*  \kron{\nvar} \mathbf{k}\S.
\end{eqnarray}
\end{subequations} 
\end{definition}

The matrices $\Lb^{\{a\}}$, $a \in \{\slow,\fast\}$, can be chosen as exact Jacobians (ROS schemes) or as arbitrary approximations to the Jacobians (ROW schemes).
Method~\eqref{eqn:MROS/MROW} is a GARK ROS/ROW scheme with step $H$ defined by the Butcher tableau \eqref{eqn:GARK-Rosenbrock-butcher}:
\begin{subequations}
\label{eqn:mrGARK-ROW-butcher}
\begin{equation}
\label{eqn:mrGARK-ROW-butcher-Alpha}
\renewcommand{\arraystretch}{1.5}
\begin{array}{c|c}
\boldsymbol{\upalpha}\FF & \boldsymbol{\upalpha}\FS \\ \hline
\boldsymbol{\upalpha}\SF &\boldsymbol{\upalpha}\SS \\ \Xhline{1.5pt}
{\b}\F  \,\!^T & {\b}\S  \,\!^T
\end{array} 
~~  \raisebox{-10pt}{$\coloneqq$} ~~ 
\raisebox{19pt}{$
\begin{array}{ccc|c}  
 \sfrac{1}{\mratio} \rosalpha\FFL[1]                      & \cdots & 0 & \rosalpha\FSL[1]  \\
\vdots                           & \ddots &  \vdots & \vdots  \\
 \sfrac{1}{\mratio}\one\F  \rosb\FL[1]*  & \ldots &  \sfrac{1}{\mratio} \rosalpha\FFL[\mratio] &\rosalpha\FSL[\mratio] \\
\hline 
 \sfrac{1}{\mratio} \rosalpha\SFL[1]  & \cdots &  \sfrac{1}{\mratio}\rosalpha\SFL[\mratio] & \rosalpha \SS   \\   \Xhline{1.5pt} 
  \sfrac{1}{\mratio} \rosb\FL[1]*  & \ldots &  \sfrac{1}{\mratio} \rosb\FL[\mratio]* & \rosb\S*  
\end{array}
$}\raisebox{-10pt}{,}
\end{equation}
\begin{equation}
\label{eqn:mrGARK-ROW-butcher-Gamma}
\renewcommand{\arraystretch}{1.5}
\begin{array}{c|c}
\boldsymbol{\upgamma}\FF & \boldsymbol{\upgamma}\FS  \\ \hline
\boldsymbol{\upgamma}\SF &\boldsymbol{\upgamma}\SS 
\end{array} 
~~  \coloneqq ~~ 
\raisebox{18pt}{$
\begin{array}{ccc|c}  
 \sfrac{1}{\mratio} \rosgamma\FFL[1]                        & \cdots & 0 & \rosgamma\FSL[1]  \\
\vdots                        & \ddots &  \vdots  & \vdots  \\
 0   & \ldots &  \sfrac{1}{\mratio} \rosgamma\FFL[\mratio] &\rosgamma\FSL[\mratio] \\
\hline 
 \sfrac{1}{\mratio} \rosgamma\SFL[1]  & \cdots &  \sfrac{1}{\mratio} \rosgamma\SFL[\mratio] & \boldsymbol\rosgamma \SS  
 \end{array}
 $},
\end{equation}
%
%
and $\boldsymbol{\upbeta}^{\{a,b\}} \coloneqq \boldsymbol{\upalpha}^{\{a,b\}} + \boldsymbol{\upgamma}^{\{a,b\}}$ for $a \in \{\slow,\fast\}$.
\end{subequations}

\begin{remark}[Intermediate micro-step solutions]
\label{rem:intermediate-fast-solutions}
Define the intermediate solutions at the end of the fast micro-steps as:
\begin{equation} 
\label{eqn:intermediate-fast-solutions}
\widetilde{\y}_{n-1+\lambda/\mratio} \coloneqq \y_{n-1}  + \sum_{\ell=1}^{\lambda} \rosb\FL[\ell]*  \kron{\nvar} \mathbf{k}\FL[\ell].
\end{equation}
From \eqref{eqn:MROS/MROW-KF} the fast stages are computed as:
\begin{eqnarray*}
\mathbf{k}\FL[\lambda] &=& h\,\fun\F\Bigl( \one_{s\F} \otimes \widetilde{\y}_{n-1+(\lambda-1)/\mratio} +
 \rosalpha\FFL[\lambda]  \kron{\nvar} \mathbf{k}\FL[\lambda]  + 
 \rosalpha\FSL[\lambda]  \kron{\nvar} \mathbf{k}\S \Bigr) \\
 && + \Bigl( \Id_{s\F \times s\F} \otimes h\,\Lb\F \Bigr) \,
\Bigl( 
\rosgamma\FFL[\lambda]  \kron{\nvar} \mathbf{k}\FL[\lambda]
+ \rosgamma\FSL[\lambda] \kron{\nvar} \mathbf{k}\S \Bigr),
\end{eqnarray*}
and from \eqref{eqn:MROS/MROW-Sol} the next step solution as:
\begin{eqnarray*}
\y_{n} &=& \widetilde{\y}_{n} +  \rosb\S*  \kron{\nvar} \mathbf{k}\S.
\end{eqnarray*}
\end{remark}

\begin{remark}[Non-uniform micro-steps]
\label{rem:microsteps}
Definition \ref{def:MGARK-ROS/ROW} can be immediately extended to accommodate nonuniform micro-steps $h_\lambda$ with $\sum_{\lambda=1}^\mratio h_\lambda = H$ by using $h_\lambda$ in each fast step \eqref{eqn:MROS/MROW-KF}, and scaling columns $\lambda$ of the Butcher tableaus \eqref{eqn:mrGARK-ROW-butcher-Alpha} and \eqref{eqn:mrGARK-ROW-butcher-Gamma} by $h_\lambda/H$ instead of  $1/\mratio$.
\end{remark}

\begin{remark}[Pure multirate approach]
\label{rem:pure-multirate}
In this paper we focus on ``pure'' multirate methods where all micro-steps use the same fast base method, i.e., 
\begin{equation}
\label{eqn:pure-multirate}
\left(  \rosb\FL[\lambda],\rosalpha\FFL[\lambda],\rosgamma\FFL[\lambda] \right) \equiv \left(  \rosb\F,\rosalpha\F,\rosgamma\F \right), \quad \lambda = 1,\dots,\mratio. 
\end{equation}
Nevertheless, the general structure of Definition \ref{def:MGARK-ROS/ROW} remains of interest as it allows to build special types of fast-slow couplings.
\end{remark}

\subsection{Multirate GARK-ROS/ROW for component-wise splitting}

Consider a two-way component partitioned system \eqref{eqn:multirate-ode} of the form
\begin{equation} 
\label{eqn:component-ode}
\renewcommand{\arraystretch}{1.5}
\y' = \begin{bmatrix} \y\F \\ \y\S \end{bmatrix}' = \begin{bmatrix} \fun\F(\y\F, \y\S) \\ \fun\S(\y\F, \y\S) \end{bmatrix} = \fun(\y), \quad  \begin{array}{ll} \y\F(t_0)&=\y\F_0 \in \Re^{\nvar\F}, \\ \y\S(t_0)&=\y\S_0 \in \Re^{\nvar\S}, \end{array}
\end{equation}
again, with a slow part $\fun\S$ that is computationally expensive and a fast part $\fun\F$ that is inexpensive to evaluate. 

\begin{definition}[Multirate GARK-ROS/ROW method for component partitioned systems] One macro-step of a MR-GARK-ROS/ROW method applied to \eqref{eqn:component-ode} with macro-step $H$ and $\mratio$ equal micro-steps $h = H/\mratio$ reads:
\begin{subequations}
\label{eqn:MROS/MROW-component}
\begin{eqnarray}
\mathbf{k}\FL[\lambda] &=& h\,\fun\F\Bigl( \one_{s\F} \otimes \y\F_{n-1}  + 
\sum_{\ell=1}^{\lambda-1}  \rosb\FL[\lambda]*  \kron{\fast}\mathbf{k}\FL[\ell] +
 \rosalpha\FFL[\lambda]  \kron{\fast} \mathbf{k}\FL[\lambda],  \nonumber \\
\label{eqn:MROS/MROW-KF-component}
 & & \qquad\qquad  \one_{s\F} \otimes \y\S_{n-1} + \rosalpha\FSL[\lambda]  \kron{\slow} \mathbf{k}\S \Bigr) \\
&& + \Bigl( \Id_{s\F \times s\F} \otimes h \Lb\FF \Bigr) \,
\Bigl( \rosgamma\FFL[\lambda]  \kron{\fast} \mathbf{k}\FL[\lambda] \Bigr)  \nonumber \\
&& + \Bigl( \Id_{s\F \times s\F} \otimes h \Lb\FS \Bigr) \,
\Bigl( \rosgamma\FSL[\lambda] \kron{\slow} \mathbf{k}\S \Bigr),  \nonumber \\[3pt]
& &\quad  \textnormal{for}\quad \lambda=1,\ldots,\mratio;  \nonumber \\
\mathbf{k}\S &=& H\,\fun\S\Bigl( \one_{s\S} \otimes \y\F_{n-1}  + 
\sum_{\lambda=1}^\mratio \rosalpha\SFL[\lambda] \kron{\fast} \mathbf{k}\FL[\lambda],  \nonumber \\
\label{eqn:MROS/MROW-KS-component}
&&\qquad \qquad \one_{s\S} \otimes \y\S_{n-1}  +    \rosalpha\SS  \kron{\slow} \mathbf{k}\S \Bigr)   \\
& &  +  \bigl( \Id_{s\S \times s\S} \otimes H \Lb\SF) \, \Bigl(
\sum_{\lambda=1}^\mratio\rosgamma\SFL[\lambda]  \kron{\fast} \mathbf{k}\FL[\lambda] \Bigr)  \nonumber \\
&& +  \bigl( \Id_{s\S \times s\S} \otimes H \Lb\SS) \, \Bigl(
 \rosgamma\SS  \kron{\slow} \mathbf{k}\S \Bigr);  \nonumber \\
\y\F_{n} &=& \y\F_{n-1} + \sum_{\lambda=1}^{\mratio} \rosb\FL[\lambda]*  \kron{\fast} \mathbf{k}\FL[\lambda]; \\
\y\S_{n} &=& \y\S_{n-1} + \rosb\S*  \kron{\slow} \mathbf{k}\S,
\end{eqnarray}
\end{subequations} 
\end{definition}
where $\kron{\fast}$ and $\kron{\slow}$ are the notation \eqref{eqn:K-matrix-notation} for $\nvar\F$ and $\nvar\S$, respectively.
In case of GARK-ROS methods we define 
\ifnum \value{book}=0
$\Lb^{\{a,b\}} \coloneqq \partial \fun^{\{a\}}/\partial \y^{\{b\}}(\y\F_{n-1}, \y\S_{n-1})$, 
$a,b \in \{ \fast, \slow \}$,
\else
\begin{equation*}
\Lb^{\{a,b\}} \coloneqq \frac{\partial \fun^{\{a\}}}{\partial \y^{\{b\}}}(\y\F_{n-1}, \y\S_{n-1}), \quad a,b \in \{ \fast, \slow \},
\end{equation*}
\fi
and in case of GARK-ROW methods these are arbitrary approximations of the Jacobian blocks.

\section{Order conditions for MR-GARK-ROS/ROW schemes}
\label{sec:order-conditions}
%
In the previous section we have seen that MR-GARK-ROS/ROW  schemes can be interpreted as partitioned GARK-ROS/ROW schemes. Thus the order conditions of these multirate schemes can be easily derived form the underlying GARK-ROS/ROW schemes derived in ~\cite{Sandu_2021_GARK-ROS}. 
In the following we focus on pure multirate methods \eqref{eqn:pure-multirate}. We consider only order conditions up to order four, as multirate schemes fit better to lower error tolerances and thus lower order schemes. Lower error tolerances are sufficient for many real-world applications, if they are in the range of the respective modeling and measurement errors of model parameters.

For exact Jacobians the order conditions for GARK-ROS schemes are \eqref{eqn:GARK-ROS-order-conditions}.
For GARK-ROS schemes with time-lagged Jacobian information, in addition to the order conditions~\eqref{eqn:GARK-ROS-order1-conditions}--\eqref{eqn:GARK-ROS-order3-conditions}, the order conditions
\begin{eqnarray}
\label{eqn:GARK-ROW-time-lagged-order3-conditions}
\textnormal{time-lagged order 3: } \quad&&
\b^{\{m\}}\,\!^T  \, \c^{\{m,n\}}  = \sfrac{1}{2}
\quad  \textnormal{for } m,n \in \{\fast,\slow\}.
\end{eqnarray}
have to be fulfilled for order three. For general approximations of the Jacobians, the order conditions of GARK-ROW schemes are given by  \eqref{eqn:GARK-ROW-order-conditions}.
%
\subsection{Internal consistency of MR-GARK-ROS/ROW schemes}
We consider methods \eqref{eqn:MROS/MROW} of pure multirate type \eqref{eqn:pure-multirate}. 
\ifnum \value{book}=1
The vectors \eqref{eqn:terminal-vectors} are:
\begin{equation}
\label{eqn:terminal-vectors-multirate}
\begin{split}
\c\SS &\coloneqq \pmb{\upalpha}\SS \,  \one = \rosc\SS, \\
\c\SF &\coloneqq \pmb{\upalpha}\SF \,  \one = \sfrac{1}{\mratio} \sum_{\lambda=1}^{\mratio} \rosc\SFL[\lambda], \\
\c\FS &\coloneqq \pmb{\upalpha}\FS \,  \one = \big[ \rosc\FSL[\lambda] \big]_{1 \le \lambda \le \mratio}, \\
\c\FF &\coloneqq \pmb{\upalpha}\FF \, \one = \Big[ \sfrac{\lambda-1}{\mratio}\, \one\F +  \sfrac{1}{\mratio}\, \rosc\FF \Big]_{1 \le \lambda \le \mratio}, \\
\g\SS &\coloneqq \pmb{\upgamma}\SS \,  \one = \rosg\SS, \\
\g\SF &\coloneqq \pmb{\upgamma}\SF \,  \one = \sfrac{1}{\mratio} \sum_{\lambda=1}^{\mratio} \rosg\SFL[\lambda], \\
\g\FS &\coloneqq \pmb{\upgamma}\FS \,  \one = \big[ \rosg\FSL[\lambda] \big]_{1 \le \lambda \le \mratio}, \\
\g\FF &\coloneqq \pmb{\upgamma}\FF \, \one = \Big[ \sfrac{1}{\mratio}\, \rosg\FF \Big]_{1 \le \lambda \le \mratio}.
\end{split}
\end{equation}
\fi
The internal consistency conditions \eqref{eqn:ROW-internal-consistency} read:
\begin{subequations}
\label{eqn:internal-consistency-multirate}
\begin{eqnarray}
\label{eqn:internal-consistency-multirate-csf}
&& \sfrac{1}{\mratio} \sum_{\lambda=1}^{\mratio} \rosc\SFL[\lambda] = \rosc\SS; \\
\label{eqn:internal-consistency-multirate-cfs}
&& \rosc\FSL[\lambda] = \sfrac{\lambda-1}{\mratio}\, \one\F +  \sfrac{1}{\mratio}\, \rosc\FF, \qquad \lambda = 1, \dots,  \mratio; \\
\label{eqn:internal-consistency-multirate-gsf}
&& \sfrac{1}{\mratio} \sum_{\lambda=1}^{\mratio} \rosg\SFL[\lambda] = \rosg\SS; \\
\label{eqn:internal-consistency-multirate-gfs}
&& \rosg\FSL[\lambda] = \sfrac{1}{\mratio}\, \rosg\FF, \qquad \lambda = 1, \dots,  \mratio.
\end{eqnarray}
\end{subequations}
From here we have that:
\begin{equation}
\label{eqn:internal-consistency-multirate-beta}
\begin{split}
& \sfrac{1}{\mratio} \sum_{\lambda=1}^{\mratio} \rose\SFL[\lambda] = \rose\SS; \\
& \rose\FSL[\lambda] = \sfrac{\lambda-1}{\mratio}\, \one\F +  \sfrac{1}{\mratio}\, \rose\FF, \qquad \lambda = 1, \dots,  \mratio.
\end{split}
\end{equation}

In what follows it is convenient to use the following notation for the coupling coefficients  averaged across all micro-steps:
\begin{equation}
\label{eqn:average-notation}
\rosalpha\FSL[\Sigma_k] \coloneqq \sum_{\lambda=1}^\mratio  (\lambda-1)^k \, \rosalpha\FSL[\lambda],
\end{equation}
with similar notations used for $\rosalpha\SFL[\cdot]$, $\rosgamma\FSL[\cdot]$, and $\rosgamma\SFL[\cdot]$.

\subsection{Order four conditions for internally consistent MR-GARK-ROS schemes}

We consider internally consistent MR-GARK-ROS schemes \eqref{eqn:internal-consistency-multirate}. It is immediate that, if the base methods are second order ROS schemes, the MR-GARK-ROS has also second order \eqref{eqn:GARK-ROW-order2-conditions}.

We turn our attention to the MR-GARK-ROS order three conditions \eqref{eqn:GARK-ROS-order3-conditions}. Assume that each of the base methods is a third order ROS scheme. The remaining order three coupling conditions are as follows:
\begin{subequations}
\label{eqn:MGARK-ROS-order3} 
\begin{eqnarray}
\label{eqn:MGARK-ROS-order3-a} 
\rosb\S* \, \left( \rosbeta\SFL[\Sigma_1]\, \one\F +  \rosbeta\SFL[\Sigma_0]\,\rose\FF  \right) &=& \sfrac{\mratio^2}{6}, \\
\label{eqn:MGARK-ROS-order3-b} 
\rosb\F* \,\rosbeta\FSL[\Sigma_0] \, \rose\SS &=& \sfrac{\mratio}{6}.
\end{eqnarray} 
\end{subequations}
%

\begin{remark}[Order three conditions for time-lagged Jacobian approach]
\label{remark.time-lagged}
Consider an internally consistent MR-GARK-ROS scheme of order three. The additional order three conditions for time-lagged Jacobians~\eqref{eqn:GARK-ROW-time-lagged-order3-conditions} reduce to:
\begin{eqnarray*}
\b\S*  \,\c\SS &=& \sfrac{1}{2}, \quad
\b\F*\, \c\FF = \sfrac{1}{2}
\quad \Leftrightarrow \quad
\rosb\S*  \,\rosc\SS = \sfrac{1}{2}, \quad
\rosb\F*\, \rosc\FF = \sfrac{1}{2}.
\end{eqnarray*}
%
%
The method maintains order three in the case of time-lagged Jacobians if and only if each base method satisfies the time-lagged Jacobian conditions  above.
%
It can be shown that equation $\rosb\F*\, \rosc\FF = 1/2$  is equivalent to the order three MR-GARK-ROS condition
$\b\F*\, (\c\FF  \times \c\FS) = 1/3$ in case of internal consistency; therefore only one additional condition
$\rosb\S*\, \rosc\SS = 1/2$ is necessary and sufficient to maintain order three when time-lagged Jacobians are used.
%
\end{remark}


We consider the MR-GARK-ROS order four conditions \eqref{eqn:GARK-ROS-order4-conditions}. Assume that each of the base methods is a fourth order ROS scheme. The remaining order four coupling conditions are as follows:
\begin{subequations}
\label{eqn:MGARK-ROS-order4}
\begin{eqnarray}
\label{eqn:MGARK-ROS-order4-a}
\rosb\S* \, \Big( ( \rosalpha\SFL[\Sigma_1]\, \one\F +  \rosalpha\SFL[\Sigma_0]\,\rose\FF )  \times \rosc\SS \Big) = \sfrac{\mratio^2}{8}, \\
\label{eqn:MGARK-ROS-order4-b}
\rosb\F* \, \Big( \rosalpha\FSL[\Sigma_1]\, \rose\SS + \big( \rosalpha\FSL[\Sigma_0]\, \rose\S\big) \times  \rosc\FF \Big) = \sfrac{\mratio^2}{8}, \\
\label{eqn:MGARK-ROS-order4-c}
\rosb\S* \, \Big( \rosbeta\SFL[\Sigma_2]\, \one\F + 2 \,\rosbeta\SFL[\Sigma_1]\, \rosc\FF + \rosbeta\SFL[\Sigma_0]\,\rosc\FF[\times 2] \Big)  = \sfrac{\mratio^3}{12}, \\ 
\label{eqn:MGARK-ROS-order4-d}
\rosb\F* \, \rosbeta\FSL[\Sigma_0]\,  \rosc\S[\times 2]  = \sfrac{\mratio}{12}, \\ 
\label{eqn:MGARK-ROS-order4-e}
\rosb\S* \,\rosbeta\SS\, \left(  \rosbeta\SFL[\Sigma_1]\, \one\F +   \rosbeta\SFL[\Sigma_0]\,\rose\FF \right) =  \sfrac{\mratio^2}{24}, \\
\label{eqn:MGARK-ROS-order4-f}
\rosb\S* \,\left( \rosbeta\SFL[\Sigma_1]\, \rose\FF +  \rosbeta\SFL[\Sigma_0]\,\rosbeta\FF\,\rose\FF \right) =  \sfrac{\mratio^3}{24}, \\
\label{eqn:MGARK-ROS-order4-g}
\rosb\S* \, \left( \sum_{\lambda=1}^\mratio \rosbeta\SFL[\lambda]\, \rosbeta\FSL[\lambda] \right) \,\rose\SS =  \sfrac{\mratio}{24}, \\
%
\label{eqn:MGARK-ROS-order4-h}
\rosb\F*  \,\rosbeta\FF\,\rosbeta\FSL[\Sigma_0] \,\rose\SS =  \sfrac{\mratio^2}{24}, \\
\label{eqn:MGARK-ROS-order4-i}
\rosb\F* \,  \rosbeta\FSL[\Sigma_0] \,\rosbeta\SS\,\rose\SS =  \sfrac{\mratio}{24}, \\
\label{eqn:MGARK-ROS-order4-j}
\rosb\F*  \,\sum_{\lambda=1}^\mratio \left( \rosbeta\FSL[\lambda] \,\rosbeta\SFL[\lambda]\, ( (\lambda-1)\, \one\F +  \rose\FF ) \right) =  \sfrac{\mratio^3}{24}.
\end{eqnarray}
\end{subequations} 

\subsection{Order three conditions for  internally consistent  MR-GARK-ROW schemes}

We consider internally consistent MR-GARK-ROW schemes \eqref{eqn:internal-consistency-multirate}. 

Assume that the base methods are second order ROW schemes. The MR-GARK-ROW second order conditions \eqref{eqn:GARK-ROW-order2-conditions} are:
\begin{eqnarray*}
\rosb\S* \, \rosc\SFL[\Sigma_0] &=& \sfrac{\mratio}{2},  \qquad
\rosb\S* \, \rosg\SFL[\Sigma_0]  = 0,  \\
\rosb\F* \, \rosc\FSL[\Sigma_0] &=& \sfrac{\mratio}{2},  \qquad
\rosb\F* \, \rosg\FSL[\Sigma_0]  = 0.
\end{eqnarray*}
These conditions are automatically satisfied for internally consistent methods.

We consider the MR-GARK-ROW order three conditions \eqref{eqn:GARK-ROW-order3-conditions}. Assume that each of the base methods is a third order ROW scheme. There are eight order three coupling conditions, as follows:
\begin{subequations}
\label{eqn:MGARK-ROW-order3}
\begin{eqnarray}
%
\rosb\S* \,  \big( \rosalpha\SFL[\Sigma_1]\, \one\F +   \rosalpha\SFL[\Sigma_0]\,\rosc\FF \big) &=& \sfrac{\mratio^2}{6}, \\
\rosb\F*\, \rosalpha\FSL[\Sigma_0]\, \rosc\SS &=& \sfrac{\mratio}{6}, \\
%
\rosb\S* \,  \big( \rosgamma\SFL[\Sigma_1]\, \one\F +   \rosgamma\SFL[\Sigma_0]\,\rosc\FF \big) &=& 0, \\
\rosb\F*\, \rosgamma\FSL[\Sigma_0] \,\rosc\SS &=& 0, \\
%
\rosb\S*  \,  \rosalpha\SFL[\Sigma_0] \, \rosg\FF &=& 0, \\
\rosb\F*\,  \rosalpha\FSL[\Sigma_0]\,\rosg\SS &=& 0, \\
%
\rosb\S*  \,  \rosgamma\SFL[\Sigma_0] \, \rosg\FF&=& 0, \\
\rosb\F*\, \rosgamma\FSL[\Sigma_0]\, \rosg\SS &=& 0. 
\end{eqnarray}
\end{subequations}

\begin{remark}
Order four coupling conditions for MR-GARK-ROW schemes can be derived in an analogous manner from the general order conditions \eqref{eqn:GARK-ROW-order4-conditions}.
\end{remark}

\section{Linear stability}
\label{sec:linear-stability}

Consider the scalar test problem
\begin{equation}
\label{eqn:scalar-dahlquist}
\y' = \lambda\S\,\y + \lambda\F\, \y.
\end{equation}
Application of the MR-GARK-ROS method \eqref{eqn:MROS/MROW} to \eqref{eqn:scalar-dahlquist} leads to the same stability equation as the application of a Multirate GARK scheme.
%
Using the notation \eqref{eqn:mrGARK-ROW-butcher} and defining 
\begin{align*}
   \B &\coloneqq \begin{bmatrix}
 \boldsymbol{\upbeta}\FF &  \boldsymbol{\upbeta}\FS \\
   \boldsymbol{\upbeta}\SF & \boldsymbol{\upbeta}\SS
   \end{bmatrix} \in \Re^{s \times s}, \quad
   \b^{T} \coloneqq \begin{bmatrix}
    \b\F* & \b\S*
   \end{bmatrix} \in \Re^{1 \times s}, \\
s & \coloneqq  s\S  + \mratio \, s\F ,  \quad
z\S  \coloneqq H\,\lambda\S, \quad 
z\F  \coloneqq H\,\lambda\F, \\
Z & \coloneqq  \operatorname{diag} \, \big\{ z\F\,\Id_{s\F \times s\F}, \ldots, 
z\F\,\Id_{s\F \times s\F}, z\S\,\Id_{s\S \times s\S}
\big\} \in \Re^{s \times s},
\end{align*}
we obtain:
\begin{equation}
\label{eqn:GARK-ROS-stability}
\begin{split}
\y_{n+1} &= R(Z)\,\y_{n}, \\
R(Z) &= 1 + \b^{T} \, \left( \Idstage - Z\,\B \right)^{-1}\,Z\,\one_{s}
 = 1 + \b^{T} \,Z\, \left( \Idstage - \B\,Z \right)^{-1}\,\one_{s},
\end{split}
\end{equation}
which is the same stability function as for a GARK scheme with tableau of coefficients $(\b,\B)$. The following definition extends immediately from GARK to MR-GARK-ROS schemes.
\begin{definition}[Stiff accuracy]
Let $\delta_{s} \in \Re^s$ be a vector with the last entry equal to one, and all other entries equal to zero.
The MR-GARK-ROS method \eqref{eqn:MROS/MROW} is called stiffly accurate if
\[
\b^{T} = \delta_{s}^{T}\,\B
\quad \Leftrightarrow \quad
{b}\S*= \delta_{s\F }^{T}\,{\beta}\FSL[\mratio] ~~\textnormal{and}~~
{b}\FL[\mratio]\,\! = \delta_{s\F }^{T}\,{\beta}\FFL[\mratio].
\]
\end{definition}
{
For a stiffly accurate MR-GARK-ROS/ROW scheme the stability function \eqref{eqn:GARK-ROS-stability} becomes:
\begin{equation}
\label{eqn:GARK-ROS-stability-SA}
\begin{split}
R(Z)  = z\F[-1]\,\mathbf{e}_{s}^{T}\,\left( Z^{-1} - \B \right)^{-1}\,\one_{s}.
\end{split}
\end{equation}
If 
    $\diag \, \big\{ \zero_{s\F}, \ldots, 
\zero_{s\F}, 1/z\S\,\Id_{s\S \times s\S}
\big\} - \B$
is nonsingular, 
then $R(Z) \to 0$ when $z\F \to \infty$.
}

For component partitioned systems we consider the following model problem \cite{Kvaerno_2000_stability-MRK}:
\begin{equation}
\label{eqn:mrgark-stability-test-matrix}
\renewcommand{\arraystretch}{1.1}
\begin{bmatrix}
\dotyf \\
\dotys
\end{bmatrix}
=
\begin{bmatrix}
\lambdaf & \etas  \\
\etaf & \lambdas
\end{bmatrix}
\,
\begin{bmatrix} \yf \\ \ys  \end{bmatrix}
\end{equation}
\ifnum \value{book}=1
The dynamics is characterized by the following coefficients:
\begin{eqnarray}
\label{eqn:mrgark-ode2d-scale-ratio}
\textnormal{scale ratio:} \quad & \mu = &\frac{|\lambda^{\{\fast\}}|}{|\lambda^{\{\slow\}}|}, \\
\label{eqn:mrgark-ode2d-coupling-coefficient}
\textnormal{coupling coefficient:} \quad & \kappa =& \frac{\etaf\, \etas}{\lambda^{\{\fast\}}\, \lambda^{\{\slow\}}}. 
\end{eqnarray}
For real coefficients the ODE~\eqref{eqn:mrgark-stability-test-matrix} is stable for $\kappa<1$.
\fi
Let $\zf=H\lambdaf$, $\zs= H \lambdas$, $\ws=H\etas$, and $\wf=H\etaf$. Application of the MGARK-ROS method~\eqref{eqn:MROS/MROW-component}, regarded as a partitioned GARK-ROS scheme according to the Butcher tableau~\eqref{eqn:mrGARK-ROW-butcher},
 advances over one step $H$ via the recurrence: 
\[
\renewcommand{\arraystretch}{1.5}
\begin{bmatrix} \yf[n] \\ \ys[n]  \end{bmatrix}
 = \mathbf{M}(\zf,\zs,\ws,\wf)\,\begin{bmatrix} \yf[n-1]  \\ \ys[n-1] \end{bmatrix},
\]
with the stability matrix: 
\begin{align*}
& \mathbf{M}(\zf,\zs,\ws,\wf) = \Id_{2 \times 2} + \begin{bmatrix} \mathbf{b}\F & \mathbf{0}_{\mratio s\F} \\ \mathbf{0}_{s\S} & \mathbf{b}\S \end{bmatrix}^T \cdot  \\
& 
\begin{bmatrix} 
\Id_{Ms\F \times Ms\F} - \zf\, \boldsymbol{\upbeta}\FF & -\ws \, \boldsymbol{\upbeta}\FS \\  
-\wf\, \boldsymbol{\upbeta}\SF                        & \Id_{s\S \times s\S} - \zs \, \boldsymbol{\upbeta}\SS  
\end{bmatrix}^{-1}\cdot  \begin{bmatrix} \zf\,\one_{\mratio s\F} & \ws\,\one_{\mratio s\F} \\ \wf\,\one_{s\S} & \zs\,\one_{s\S} \end{bmatrix}.
\end{align*}
%
One immediately sees that for one-sided coupled problems with $\ws=0$ or $\wf=0$ the stability of the base schemes guarantees the stability of the multirate schemes. 
\ifnum \value{book}=1
This is because we have
\begin{align*}
    \mathbf{M}(\zf,\zs,0,\wf)  = 
    \begin{bmatrix}
     \mathbf{R}\F & 0 \\
    \star & \mathbf{R}\S
    \end{bmatrix}, \quad
    \mathbf{M}(\zf,\zs,\ws,0)  = 
    \begin{bmatrix}
     \mathbf{R}\F & \star \\
    0 & \mathbf{R}\S
    \end{bmatrix},
\end{align*}
with $\mathbf{R}\F$, $\mathbf{R}\S$ the stability functions of the base fast and slow schemes, respectively:
\begin{align*}
  \mathbf{R}\F & = 1+ {\mathbf{b}\F}^\top \left(\Id_{Ms\F \times Ms\F} - \zf\, \boldsymbol{\upbeta}\FF \right)^{-1} 
     \zf\,\one_{\mratio s\F}, \\
   \mathbf{R}\S & =  1+ {\mathbf{b}\S}^\top \left( \Id_{s\S \times s\S} - \zs \, \boldsymbol{\upbeta}\SS \right)^{-1} \zs\,\one_{s\S}.
\end{align*}
\fi

\section{Coupling the fast and slow systems in a computationally-efficient manner}
\label{sec:coupling-strategies}

In traditional Rosenbrock methods the coefficient matrix $\rosalpha$ is strictly lower triangular, and the matrix $\rosgamma$  lower triangular with equal diagonal entries. Due to this structure the stages $\mathbf{k}_i$ are evaluated sequentially in a decoupled manner, each stage computation is only implicit in the current stage.

Multirate GARK ROS/ROW schemes compute both slow and fast stage vectors. We call a stage computation ``decoupled'' if it is implicit in only the current stage $\mathbf{k}\S_i$ or $\mathbf{k}\FL[\lambda]_i$. We call computations ``coupled'' if one (or more) slow stages, and one (or more) fast stages are computed together by solving a single large system of linear equations. Computational efficiency of multirate methods relies on evaluating less frequently the expensive slow part. Consequently, an efficient multirate method keeps the coupling at a minimum.

In this paper we construct multirate GARK ROS/ROW schemes where the base methods are Rosenbrock(-W) schemes with matrices $\rosalpha\FFL[\lambda]$, $\rosalpha\SS$ strictly lower triangular, and matrices  $\rosgamma\FFL[\lambda]$ and $\rosgamma\SS$  lower triangular. 
From the Butcher tableau \eqref{eqn:mrGARK-ROW-butcher} compute the coupling structure matrix
\begin{equation}
\label{eqn:stucture-matrix}
\mathbf{S} \coloneqq \Bigl( |\pmb{\upalpha}\SF| + |\pmb{\upgamma}\SF| \Bigr)^T \times \Bigl( |\pmb{\upalpha}\FS| + |\pmb{\upgamma}\FS| \Bigr) \in \Re^{\mratio s\F \times s\S},
\end{equation}
where $| \cdots |$ takes element-wise absolute values, and $\times$ is the element-wise product.
We make the following observations:
\begin{itemize}

\item In order to compute stages sequentially,  in a completely decoupled manner, it must hold that $\mathbf{S} = \zero$.

\item The non-zero entries in this matrix $\mathbf{S} $ correspond to slow and fast stages that are computed together, in a coupled manner. Specifically, if element in row $(\lambda,i)$ and column  $j$ is non-zero then stages $\mathbf{k}\FL[\lambda]_i$ and $\mathbf{k}\S_j$ are computed by solving a joint linear system.

\item Note that internal consistency equation \eqref{eqn:ROW-internal-consistency} for $\g$ requires that at least one slow and one fast stage are computed together in a coupled manner.

\end{itemize}

\begin{example}[Second order, two-rate method]
\label{example:second-order-M2}
Consider the following example using $\mratio=2$ and two-stage base methods:
\begin{equation*}
\renewcommand{\arraystretch}{1.5}
\scalebox{0.6}{$
\begin{array}{cc:cc|cc}  
0     &     0 &     0  & 0 & 0 & 0 \\
{\scriptstyle \frac{1}{2}}\, \alpha\FFL[1]_{2,1}     &     0 &     0     & 0 & \alpha\FSL[1]_{2,1} & \alpha\FSL[1]_{2,2}  \\
\hdashline
{\scriptstyle \frac{1}{2}}\, b\FL[1]_1 & {\scriptstyle \frac{1}{2}}\,b\FL[1]_2 & 0  &0 & \alpha\FSL[2]_{1,1} & \alpha\FSL[2]_{1,2}  \\
{\scriptstyle \frac{1}{2}}\,b\FL[1]_1 & {\scriptstyle \frac{1}{2}}\,b\FL[1]_2 & {\scriptstyle \frac{1}{2}}\,\alpha\FFL[2]_{2,1}     &0 & \alpha\FSL[2]_{2,1} & \alpha\FSL[2]_{2,2}  \\
\hline
0 & 0 & 0 & 0 & 0 & 0 \\  
{\scriptstyle \frac{1}{2}}\,\alpha\SFL[1]_{2,1} & {\scriptstyle \frac{1}{2}}\,\alpha\SFL[1]_{2,2} & {\scriptstyle \frac{1}{2}}\,\alpha\SFL[2]_{2,1}  & {\scriptstyle \frac{1}{2}}\,\alpha\SFL[2]_{2,2} & \alpha\SS_{2,1}   & 0 \\   
\Xhline{1.5pt}
{\scriptstyle \frac{1}{2}}\,b\FL[1]_1 & {\scriptstyle \frac{1}{2}}\,b\FL[1]_2 & {\scriptstyle \frac{1}{2}}\,b\FL[2]_1 & {\scriptstyle \frac{1}{2}}\,b\FL[2]_2 & b\S_{1} & b\S_{2} \\  
\end{array}
$}
\raisebox{-4pt}{$\,,$}
\quad
\raisebox{6pt}{
\scalebox{0.6}{$
\renewcommand{\arraystretch}{1.5}
\begin{array}{cc:cc|cc}  
{\scriptstyle \frac{1}{2}}\, \gamma\FFL[1]   &     0 &     0  & 0 & \gamma\FSL[1]_{1,1} & \gamma\FSL[1]_{1,2}  \\
{\scriptstyle \frac{1}{2}}\, \gamma\FFL[1]_{2,1}     &  {\scriptstyle \frac{1}{2}}\, \gamma\FFL[1]  &     0     & 0 & \gamma\FSL[1]_{2,1} & \gamma\FSL[1]_{2,2}  \\
\hdashline
0 & 0 & {\scriptstyle \frac{1}{2}}\,\gamma\FFL[2]    &0 & \gamma\FSL[2]_{1,1} & \gamma\FSL[2]_{1,2}  \\
0 & 0 & {\scriptstyle \frac{1}{2}}\,\gamma\FFL[2]_{2,1}     & {\scriptstyle \frac{1}{2}}\,\gamma\FFL[2]   & \gamma\FSL[2]_{2,1} & \gamma\FSL[2]_{2,2}  \\
\hline
{\scriptstyle \frac{1}{2}}\,\gamma\SFL[1]_{1,1} & {\scriptstyle \frac{1}{2}}\,\gamma\SFL[1]_{1,2} & {\scriptstyle \frac{1}{2}}\,\gamma\SFL[2]_{1,1}  & {\scriptstyle \frac{1}{2}}\,\gamma\SFL[2]_{1,2} & \gamma\SS   & 0 \\  
{\scriptstyle \frac{1}{2}}\,\gamma\SFL[1]_{2,1} & {\scriptstyle \frac{1}{2}}\,\gamma\SFL[1]_{2,2} & {\scriptstyle \frac{1}{2}}\,\gamma\SFL[2]_{2,1}  & {\scriptstyle \frac{1}{2}}\,\gamma\SFL[2]_{2,2} & \gamma\SS_{2,1}   & \gamma\SS  
\end{array}
$}}
\raisebox{-4pt}{$\,.$}
\end{equation*}
We conveniently choose $\alpha\SFL[\lambda]_{1,j}=\alpha\FSL[\lambda]_{1,j}=0$ for all $j$ and $\lambda$ such as to satisfy the first internal consistency conditions. The coupling structure matrix \eqref{eqn:stucture-matrix} is:
\begin{equation}
\label{eqn:stucture-matrix-example}
\mathbf{S} =
\scalebox{0.75}{$
\sfrac{1}{2}\, \begin{bmatrix}
|\gamma\SFL[1]_{1,1}|\cdot|\gamma\FSL[1]_{1,1}| & (|\alpha\SFL[1]_{2,1}|  + |\gamma\SFL[1]_{2,1}|)\cdot|\gamma\FSL[1]_{1,2}| \\
|\gamma\SFL[1]_{1,2}|\cdot (|\alpha\FSL[1]_{2,1}|+|\gamma\FSL[1]_{2,1}|) & (|\alpha\SFL[1]_{2,2}| + |\gamma\SFL[1]_{2,2}|) \cdot (|\alpha\FSL[1]_{2,2}|+|\gamma\FSL[1]_{2,2}|) \\
|\gamma\SFL[2]_{1,1}|\cdot (|\alpha\FSL[2]_{1,1}|+|\gamma\FSL[2]_{1,1}|)  & (|\gamma\SFL[2]_{2,1}|  + |\alpha\SFL[2]_{2,1}|)\cdot (|\alpha\FSL[2]_{1,2}|+|\gamma\FSL[2]_{1,2}|) \\
|\gamma\SFL[2]_{1,2}|\cdot (|\alpha\FSL[2]_{2,1}|+|\gamma\FSL[2]_{2,1}|) & (|\gamma\SFL[2]_{2,2}|   + |\alpha\SFL[2]_{2,2}|)\cdot (|\alpha\FSL[2]_{2,2}|+|\gamma\FSL[2]_{2,2}|)
\end{bmatrix}.
$}
\end{equation}

\end{example}

\subsection{IMEX approach}
\label{subsec:coupling-imex}
If one chooses $\rosgamma\FFL[\lambda] = \zero$ and  $\rosgamma\FSL[\lambda] = \zero$ then the fast component is integrated with a Runge-Kutta method; this method is explicit if $\rosalpha\FFL[\lambda]$ are strictly lower triangular matrices. For a decoupled computation one needs to select the coupling coefficients $\rosalpha\FSL[\lambda]$, $\rosalpha\SFL[\lambda]$, and $\rosgamma\SFL[\lambda]$ such that the matrix $\mathbf{S} = \zero$.
Using notation \eqref{eqn:intermediate-fast-solutions}, the fast stages are computed as:
\begin{eqnarray*}
\mathbf{k}_i\FL[\lambda] &=& h \,\fun\F\Bigl( \widetilde{\y}_{n-1+(\lambda-1)/\mratio} + 
\sum_{j=1}^{s\F} \alpha\FFL[\lambda]_{i,j}  \, \mathbf{k}_j\FL[\lambda]  + 
\sum_{j=1}^{s\S} \alpha\FSL[\lambda]_{i,j}  \, \mathbf{k}\S_j  \Bigr).
\end{eqnarray*}

The corresponding Butcher tableau \eqref{eqn:mrGARK-ROW-butcher} is:
\begin{equation}
\label{eqn:mrGARK-ROW-butcher-IMEX}
\renewcommand{\arraystretch}{2.0}
\scalebox{0.7}{$
\begin{array}{cccc|cccc}  
 \sfrac{1}{\mratio} \rosalpha\FFL[1]      &          0                   & \cdots & 0 & \rosalpha\FSL[1]  \\
 \sfrac{1}{\mratio} \one  \b\FL[1]* &  \sfrac{1}{\mratio} \rosalpha\FFL[2]        & \cdots & 0 &  \rosalpha\FSL[2]  \\
\vdots                     &                             & \ddots &   & \vdots  \\
 \sfrac{1}{\mratio}\one  \b\FL[1]* &  \sfrac{1}{\mratio}\one  \b\FL[2]*   & \ldots &  \sfrac{1}{\mratio} \rosalpha\FFL[\mratio] &\rosalpha\FSL[\mratio] \\
\hline 
 \sfrac{1}{\mratio} \rosalpha\SFL[1] &  \sfrac{1}{\mratio}\rosalpha\SFL[2] & \cdots &  \sfrac{1}{\mratio}\rosalpha\SFL[\mratio] & \rosalpha \SS   
\end{array}
$}\raisebox{-25pt}{\,,}
\quad
\scalebox{0.7}{$
\begin{array}{cccc|cccc}  
0    &          0                   & \cdots & 0 & 0 \\
0 & 0  & \cdots & 0 & 0 \\
\vdots    & \vdots & \ddots &  \vdots & \vdots  \\
 0 & 0   & \ldots & 0 &0 \\
\hline 
 \sfrac{1}{\mratio} \rosgamma\SFL[1] &  \sfrac{1}{\mratio} \rosgamma\SFL[2] & \cdots &  \sfrac{1}{\mratio} \rosgamma\SFL[\mratio] & \rosgamma \SS  
 \end{array}
 $}
 \raisebox{-25pt}{\,.}
\end{equation}

\subsection{Compound-first-step approach: coupling the  macro-step with the first micro-step}
\label{sec:coupled-first-step}
%
We consider base methods with the same number of stages 
$s\F=s\S = s$.
Moreover, we set the coupling coefficients $\rosalpha\SFL[\lambda] = \rosgamma\SFL[\lambda] = 0$ for $\lambda = 2, \dots,\mratio$. 
 The resulting Butcher tableau \eqref{eqn:mrGARK-ROW-butcher} is:
\begin{equation}
\label{eqn:mrGARK-ROW-butcher-compound-step}
\renewcommand{\arraystretch}{2.0}
\scalebox{0.7}{$
\begin{array}{cccc|cccc}  
 \sfrac{1}{\mratio} \rosalpha\FFL[1]      &          0                   & \cdots & 0 & \rosalpha\FSL[1]  \\
 \sfrac{1}{\mratio} \one  \b\FL[1]* &  \sfrac{1}{\mratio} \rosalpha\FFL[2]        & \cdots & 0 &  \rosalpha\FSL[2]  \\
\vdots    & \vdots & \ddots &  \vdots & \vdots  \\
 \sfrac{1}{\mratio}\one  \b\FL[1]* &  \sfrac{1}{\mratio}\one  \b\FL[2]*   & \ldots &  \sfrac{1}{\mratio} \rosalpha\FFL[\mratio] &\rosalpha\FSL[\mratio] \\
\hline 
 \sfrac{1}{\mratio} \rosalpha\SFL[1] & 0 & \cdots & 0 & \rosalpha \SS   
\end{array}
$}\raisebox{-25pt}{\,,}
\quad
\scalebox{0.7}{$
\begin{array}{cccc|cccc}  
 \sfrac{1}{\mratio} \rosgamma\FFL[1]      &          0                   & \cdots & 0 & \rosgamma\FSL[1]  \\
0 &  \sfrac{1}{\mratio} \rosgamma\FFL[2]        & \cdots & 0 &  \rosgamma\FSL[2]  \\
\vdots    & \vdots & \ddots &  \vdots & \vdots  \\
 0 & 0   & \ldots &  \sfrac{1}{\mratio} \rosgamma\FFL[\mratio] &\rosgamma\FSL[\mratio] \\
\hline 
 \sfrac{1}{\mratio} \rosgamma\SFL[1] & 0 & \cdots & 0 & \rosgamma \SS  
 \end{array}
 $}
 \raisebox{-25pt}{\,.}
 \end{equation}
 The coefficient matrices $\alpha\FFL[1]$, $\alpha\FSL[1]$, $\alpha\SFL[1]$, $\alpha\SS$ are chosen strictly lower triangular.
The coefficient matrices $\gamma\FFL[1]$, $\gamma\FSL[1]$, $\gamma\SFL[1]$, and $\gamma\SS$ are chosen lower triangular, with equal diagonal entries: $\gamma_{i,i}\SS \coloneqq \gamma\SS, \gamma_{i,i}\FFL[1] \coloneqq \gamma\FFL[1], \gamma_{i,i}\SFL[1] \coloneqq \gamma\SFL[1], \gamma_{i,i}\FSL[1] \coloneqq \gamma\FSL[1]$ for $i=1,\dots,s$.

 In the structure matrix \eqref{eqn:stucture-matrix} the entries $\mathbf{S}_{i,i} \ne 0$ for $i=1,\dots, s\S$. This means that each slow stage $\mathbf{k}\S_i$ and the corresponding fast stage of the first micro-step $\mathbf{k}\FL[1]_i$ are computed in a coupled manner, by solving a full coupled system of linear equations. 
\ifnum \value{book}=1 
Assume that the first  stages of the first fast micro-step $\mathbf{k}_1\FL[1] \dots \mathbf{k}_{i-1}\FL[1]$, and the first slow stages $\mathbf{k}_1\S \dots \mathbf{k}_{i-1}\S$, have been computed. Stages $i$ are computed together in a coupled manner, as follows:
\begin{eqnarray*}
&&\begin{bmatrix} \Id - h\, \gamma\FFL[1]_{i,i}\, \Lb\F & -  h\, \Lb\F \, \gamma\FSL[1]_{i,i} \\
- H \, \gamma\SFL[1]_{i,i}\, \Lb\S & \Id -  H\,\gamma\SS_{i,i}\, \Lb\S \,\end{bmatrix} \, \begin{bmatrix} \mathbf{k}_i\FL[1] \\ \mathbf{k}\S_i \end{bmatrix} \\
&&= \begin{bmatrix}  h \,\fun\F\Bigl( \y_{n-1} + 
\sum_{j=1}^{i-1} \alpha\FFL[1]_{i,j}  \, \mathbf{k}_j\FL[1]  + 
\sum_{j=1}^{i-1} \alpha\FSL[1]_{i,j}  \, \mathbf{k}\S_j  \Bigr) \\
\qquad + h\, \Lb\F \, \Bigl( 
\sum_{j=1}^{i} \gamma\FFL[1]_{i,j}  \,\mathbf{k}_j \FL[1]  + 
\sum_{j=1}^{i} \gamma\FSL[1]_{i,j}  \, \mathbf{k}\S_j \Bigr)  \\[3pt]
H\,\fun\S\Bigl( \y_{n-1} +
\sum_{j=1}^{i-1} \alpha\SFL[1]_{i,j}  \, \mathbf{k}_j\FL[1]  + 
\sum_{j=1}^{i-1} \alpha\SS_{i,j}  \, \mathbf{k}\S_j  \Bigr)  \\
\qquad  +  H\, \Lb\S \, \Bigl( \sum_{j=1}^{i-1}
\gamma\SFL[1]_{i,j}\, \mathbf{k}\FL[1]_j
+ \sum_{j=1}^{i-1}\gamma\SS_{i,j} \mathbf{k}\S_j  \Bigr)
\end{bmatrix}.
\end{eqnarray*}
\fi 
Since all diagonal entries are equal to each other, only one LU decomposition of the compound matrix is necessary for computing all stage vectors $\mathbf{k}\S_i$  and $\mathbf{k}\FL[1]_i$ for $i=1,\dots,s$.

Since all slow stages are known after the first micro-step,
$\rosalpha\FSL[\lambda]$ and $\rosgamma\FSL[\lambda]$ can be full matrices for $\lambda=2,\ldots,\mratio$. For all remaining micro steps a single additional LU decomposition is necessary if $\gamma_{i,i}\FFL[\lambda]\coloneqq \gamma\FF$ is constant for all $i=1,\ldots, s\F $ and $\lambda=2,\ldots,\mratio$. 

A simple choice of coefficients for compound-first-step coupling is $\alpha\SFL[1]=\mratio\,\alpha\SS$, $\gamma\SFL[1]=\mratio\,\gamma\SS$, $\alpha\FSL[1]=(1/\mratio)\,\alpha\FFL[1]$, $\gamma\FSL[1]=(1/\mratio)\,\gamma\FFL[1]$.

\begin{example}
Consider the scheme of Example \ref{example:second-order-M2} with the following coefficients:
\begin{equation*}
\renewcommand{\arraystretch}{1.5}
\scalebox{0.6}{$
\begin{array}{cc:cc|cc}  
0     &     0 &     0  & 0 & 0 & 0 \\
{\scriptstyle \frac{1}{2}}\, \alpha\FFL[1]_{2,1}     &     0 &     0     & 0 & \alpha\FSL[1]_{2,1} & 0 \\
\hdashline
{\scriptstyle \frac{1}{2}}\, b\FL[1]_1 & {\scriptstyle \frac{1}{2}}\,b\FL[1]_2 & 0  &0 & \alpha\FSL[2]_{1,1} & \alpha\FSL[2]_{1,2}  \\
{\scriptstyle \frac{1}{2}}\,b\FL[1]_1 & {\scriptstyle \frac{1}{2}}\,b\FL[1]_2 & {\scriptstyle \frac{1}{2}}\,\alpha\FFL[2]_{2,1}     &0 & \alpha\FSL[2]_{2,1} & \alpha\FSL[2]_{2,2}  \\
\hline
0 & 0 & 0 & 0 & 0 & 0 \\  
{\scriptstyle \frac{1}{2}}\,\alpha\SFL[1]_{2,1} & 0 & 0  & 0 & \alpha\SS_{2,1}   & 0  
\end{array}
$}
\raisebox{-10pt}{$\,,$}
\quad
\scalebox{0.6}{$
\renewcommand{\arraystretch}{1.5}
\begin{array}{cc:cc|cc}  
{\scriptstyle \frac{1}{2}}\, \gamma\FFL[1]   &     0 &     0  & 0 & \gamma\FSL[1] & 0  \\
{\scriptstyle \frac{1}{2}}\, \gamma\FFL[1]_{2,1}     &  {\scriptstyle \frac{1}{2}}\, \gamma\FFL[1]  &     0     & 0 & \gamma\FSL[1]_{2,1} & \gamma\FSL[1]  \\
\hdashline
0 & 0 & {\scriptstyle \frac{1}{2}}\,\gamma\FFL[2]    &0 & \gamma\FSL[2]_{1,1} & \gamma\FSL[2]_{1,2}  \\
0 & 0 & {\scriptstyle \frac{1}{2}}\,\gamma\FFL[2]_{2,1}     & {\scriptstyle \frac{1}{2}}\,\gamma\FFL[2]   & \gamma\FSL[2]_{2,1} & \gamma\FSL[2]_{2,2}  \\
\hline
{\scriptstyle \frac{1}{2}}\,\gamma\SFL[1] & 0 & 0 & 0 & \gamma\SS   & 0 \\  
{\scriptstyle \frac{1}{2}}\,\gamma\SFL[1]_{2,1} & {\scriptstyle \frac{1}{2}}\,\gamma\SFL[1] & 0  & 0 & \gamma\SS_{2,1}   & \gamma\SS  
\end{array}
$}
\raisebox{-10pt}{$\,.$}
\end{equation*}
The coupling matrix \eqref{eqn:stucture-matrix-example} reads:
\begin{equation*}
\mathbf{S} =
\scalebox{0.75}{$
\sfrac{1}{2}\, \begin{bmatrix}
|\gamma\SFL[1]|\cdot|\gamma\FSL[1]| & 0 \\
0 &  |\gamma\SFL[1]| \cdot |\gamma\FSL[1]| \\
0 & 0 \\
0 & 0
\end{bmatrix},
$}
\end{equation*}
which indicates that $\mathbf{k}\FL[1]_1$ and $\mathbf{k}\S_1$ are computed together, and so are $\mathbf{k}\FL[1]_2$ and $\mathbf{k}\S_2$.
\ifnum \value{book}=1 
We solve the first fast stage of the first micro-step together with the first slow stage:
\begin{eqnarray*}
&& \begin{bmatrix} \Id - h\, \gamma\FFL[1]\, \Lb\F & -  h\, \Lb\F \, \gamma\FSL[1] \\
- H \, \gamma\SFL[1]\, \Lb\S & \Id -  H\,\gamma\SS\, \Lb\S \,\end{bmatrix} \, \begin{bmatrix} \mathbf{k}_1\FL[1] \\ \mathbf{k}\S_1 \end{bmatrix}  =  \begin{bmatrix} h \,\fun\F\Bigl( \y_{n-1} \Bigr) \\ H\,\fun\S\left( \y_{n-1} \right) \end{bmatrix}.
%
\end{eqnarray*}
Moreover, if all coefficients are equal $\gamma\SFL[1]_{i,i} = \gamma\FSL[1]_{i,i} = \gamma\FFL[1]  = \gamma\SS = \gamma$ the linear system to be solved has the form:
\begin{equation*}
\begin{split}
\left(\Id - h\,\gamma\,\Lb\F - H\,\gamma\,\Lb\S \right)\,(\mathbf{k}\FL[1]_1 +  \mathbf{k}\S_1) &=  \,\fun\F\left( \y_{n-1}  \right) 
+ H\,\fun\S\left( \y_{n-1} \right),
\end{split}
\end{equation*}
and then we recover individual stages from their sum using:
\begin{equation*}
\begin{split}
\mathbf{k}_1\FL[1] &= h \,\fun\F\left( \y_{n-1}  \right) + h\,\gamma\, \Lb\F \, \Bigl( 
\mathbf{k}_1 \FL[1]  + \mathbf{k}\S_1 \Bigr),  \\
\mathbf{k}\S_1 &= H\,\fun\S\left( \y_{n-1} \right)  +  H\,\gamma\, \Lb\S \, \left( \mathbf{k}\FL[1]_1 + \mathbf{k}\S_1  \right).
\end{split}
\end{equation*}
\fi 
\end{example}

\begin{remark}
The multirate ROW schemes introduced by Bartel and G\"unther~\cite {Bartel_2002_MR-W} fall into the class of multirate GARK-ROW schemes. They consider the case of time-lagged Jacobians (which differ by a term of magnitude ${\cal{O}}(H)$ from the exact Jacobian). In addition, the same order $p$ within the micro steps is demanded.
\end{remark}

\subsection{Coupling only the first fast and the first slow stage computations}
\label{sec:coupled-first-stage}
The smallest amount of coupling that allows the construction of internally consistent implicit schemes is a lighter version of the strategy discussed in section \ref{sec:coupled-first-step}, where only the first fast stage $\mathbf{k}_1\FL[1]$ and the first slow stage $\mathbf{k}\S_1$ are computed together. It is possible to select coefficients such that all subsequent stage computations are implicit in either fast or slow stages, and are computed in a decoupled manner.

\begin{example}
Consider the scheme from Example \ref{example:second-order-M2} with the following coefficients:
\begin{equation*}
\renewcommand{\arraystretch}{1.5}
\scalebox{0.6}{$
\begin{array}{cc:cc|cc}  
0     &     0 &     0  & 0 & 0 & 0 \\
{\scriptstyle \frac{1}{2}}\, \alpha\FFL[1]_{2,1}     &     0 &     0     & 0 & \alpha\FSL[1]_{2,1} & 0 \\
\hdashline
{\scriptstyle \frac{1}{2}}\, b\FL[1]_1 & {\scriptstyle \frac{1}{2}}\,b\FL[1]_2 & 0  &0 & \alpha\FSL[2]_{1,1} & 0  \\
{\scriptstyle \frac{1}{2}}\,b\FL[1]_1 & {\scriptstyle \frac{1}{2}}\,b\FL[1]_2 & {\scriptstyle \frac{1}{2}}\,\alpha\FFL[2]_{2,1}     &0 & \alpha\FSL[2]_{2,1} & \alpha\FSL[2]_{2,2}  \\
\hline
0 & 0 & 0 & 0 & 0 & 0 \\  
{\scriptstyle \frac{1}{2}}\,\alpha\SFL[1]_{2,1} & {\scriptstyle \frac{1}{2}}\,\alpha\SFL[1]_{2,2} & {\scriptstyle \frac{1}{2}}\,\alpha\SFL[2]_{2,1}  & 0 & \alpha\SS_{2,1}   & 0 
\end{array}
$}
\raisebox{-10pt}{$\,,$}
\quad
\scalebox{0.6}{$
\renewcommand{\arraystretch}{1.5}
\begin{array}{cc:cc|cc}  
{\scriptstyle \frac{1}{2}}\, \gamma\FFL[1]   &     0 &     0  & 0 & \gamma\FSL[1]_{1,1} & 0  \\
{\scriptstyle \frac{1}{2}}\, \gamma\FFL[1]_{2,1}     &  {\scriptstyle \frac{1}{2}}\, \gamma\FFL[1]  &     0     & 0 & \gamma\FSL[1]_{2,1} & 0  \\
\hdashline
0 & 0 & {\scriptstyle \frac{1}{2}}\,\gamma\FFL[2]    &0 & \gamma\FSL[2]_{1,1} & 0  \\
0 & 0 & {\scriptstyle \frac{1}{2}}\,\gamma\FFL[2]_{2,1}     & {\scriptstyle \frac{1}{2}}\,\gamma\FFL[2]   & \gamma\FSL[2]_{2,1} & \gamma\FSL[2]_{2,2}  \\
\hline
{\scriptstyle \frac{1}{2}}\,\gamma\SFL[1]_{1,1} & 0 & 0 & 0 & \gamma\SS   & 0 \\  
{\scriptstyle \frac{1}{2}}\,\gamma\SFL[1]_{2,1} & {\scriptstyle \frac{1}{2}}\,\gamma\SFL[1]_{2,2} & {\scriptstyle \frac{1}{2}}\,\gamma\SFL[2]_{2,1}  & 0 & \gamma\SS_{2,1}   & \gamma\SS  
\end{array}
$}
\raisebox{-10pt}{$\,.$}
\end{equation*}
The coupling matrix \eqref{eqn:stucture-matrix-example} has single non-zero element, ${\scriptstyle \frac{1}{2}}\,|\gamma\SFL[1]_{1,1}|\, |\gamma\FSL[1]_{1,1}|$, corresponding to first computing stages $\mathbf{k}_1\FL[1]$ and $\mathbf{k}\S_1$ in a coupled manner. Next, $\mathbf{k}_2\FL[1]$ and $\mathbf{k}_1\FL[2]$ are computed in a decoupled manner, since they only depend on the known slow stage $\mathbf{k}\S_1$. After this, $\mathbf{k}\S_2$ is computed in a decoupled manner as it does not depend on the (yet unknown) last fast stage $\mathbf{k}_2\FL[2]$.  Finally, $\mathbf{k}_2\FL[2]$ is evaluated using both slow stages.
\end{example}

\subsection{Fully decoupled approach}
%
In the completely decoupled approach each stage follows a regular Rosenbrock computation, implicit in either the fast or the slow stages, but {\it not} in both at the same time. In this case the second internal consistency conditions $\gamma\FFL[1] = \sum_j \gamma\FSL[1]_{1,j}$ do not hold unless this first stage is explicit. Consequently, the coupling order conditions for the entire method become more complex, but such methods are possible to construct.

\begin{example}
Consider the scheme from Example \ref{example:second-order-M2} with the following coefficients:
\begin{equation*}
\renewcommand{\arraystretch}{1.5}
\scalebox{0.6}{$
\begin{array}{cc:cc|cc}  
0     &     0 &     0  & 0 & 0 & 0 \\
{\scriptstyle \frac{1}{2}}\, \alpha\FFL[1]_{2,1}     &     0 &     0     & 0 & \alpha\FSL[1]_{2,1} & 0 \\
\hdashline
{\scriptstyle \frac{1}{2}}\, b\FL[1]_1 & {\scriptstyle \frac{1}{2}}\,b\FL[1]_2 & 0  &0 & \alpha\FSL[2]_{1,1} & 0  \\
{\scriptstyle \frac{1}{2}}\,b\FL[1]_1 & {\scriptstyle \frac{1}{2}}\,b\FL[1]_2 & {\scriptstyle \frac{1}{2}}\,\alpha\FFL[2]_{2,1}     &0 & \alpha\FSL[2]_{2,1} & \alpha\FSL[2]_{2,2}  \\
\hline
{\scriptstyle \frac{1}{2}}\,\alpha\SFL[1]_{1,1}  & 0 & 0 & 0 & 0 & 0 \\  
{\scriptstyle \frac{1}{2}}\,\alpha\SFL[1]_{2,1} & {\scriptstyle \frac{1}{2}}\,\alpha\SFL[1]_{2,2} & {\scriptstyle \frac{1}{2}}\,\alpha\SFL[2]_{2,1}  & 0 & \alpha\SS_{2,1}   & 0 
\end{array}
$}
\raisebox{-10pt}{$\,,$}
\quad
\scalebox{0.6}{$
\renewcommand{\arraystretch}{1.5}
\begin{array}{cc:cc|cc}  
{\scriptstyle \frac{1}{2}}\, \gamma\FFL[1]   &     0 &     0  & 0 & 0 & 0  \\
{\scriptstyle \frac{1}{2}}\, \gamma\FFL[1]_{2,1}     &  {\scriptstyle \frac{1}{2}}\, \gamma\FFL[1]  &     0     & 0 & \gamma\FSL[1]_{2,1} & 0  \\
\hdashline
0 & 0 & {\scriptstyle \frac{1}{2}}\,\gamma\FFL[2]    &0 & \gamma\FSL[2]_{1,1} & 0  \\
0 & 0 & {\scriptstyle \frac{1}{2}}\,\gamma\FFL[2]_{2,1}     & {\scriptstyle \frac{1}{2}}\,\gamma\FFL[2]   & \gamma\FSL[2]_{2,1} & \gamma\FSL[2]_{2,2}  \\
\hline
{\scriptstyle \frac{1}{2}}\,\gamma\SFL[1]_{1,1} & 0 & 0 & 0 & \gamma\SS   & 0 \\  
{\scriptstyle \frac{1}{2}}\,\gamma\SFL[1]_{2,1} & {\scriptstyle \frac{1}{2}}\,\gamma\SFL[1]_{2,2} & {\scriptstyle \frac{1}{2}}\,\gamma\SFL[2]_{2,1}  & 0 & \gamma\SS_{2,1}   & \gamma\SS  
\end{array}
$}
\raisebox{-10pt}{$\,.$}
\end{equation*}
Note the complementary sparsity structure of the off-diagonal coupling blocks. The first fast stage is that of a classical Rosenbrock method:
\begin{eqnarray*}
\Big( \Id -  h\,\gamma\FFL[1]  \, \Lb\F \Big)\,\mathbf{k}_1\FL[1] &=& h \,\fun\F\bigl( \y_{n-1} \bigr).
\end{eqnarray*}
Similarly, the first slow stage is computed in a decoupled manner:
\begin{eqnarray*}
\Big( \Id -  H\,\gamma\SS  \, \Lb\S \Big)\,\mathbf{k}_1\S &=& H\,\fun\S\Bigl(\y_{n-1}  + 
\rosalpha\SFL[1]_{1,1} \, \mathbf{k}_1\FL[1]\Bigr)   
 +  H\,\gamma\SFL[1]_{1,1}\, \Lb\S \,\mathbf{k}\FL[1]_1,
\end{eqnarray*}
and the decoupled computations continue alternating fast and slow stages.
\end{example}

\subsection{Step-predictor-corrector approach}
\label{subsec:step-predictor-corrector}
%
This approach starts with a ``predictor'' step where the slow Rosenbrock method is applied with step size $H$ to the entire system, in a classical fashion. The slow components are sufficiently accurate, but the fast components are not; for this reason we keep only the computed $\mathbf{k}\S$, but discard $\mathbf{k}\F$. The ``corrector'' re-computes $\mathbf{k}\FL[\lambda]$ for all sub-steps $\lambda$, with the small steps sizes $h$, and uses these values to construct the final solution. The Butcher tableaus \eqref{eqn:mrGARK-ROW-butcher} read:
\begin{equation}
\label{eqn:MROS/MROW-SPC-butcher}
\renewcommand{\arraystretch}{2.0}
\scalebox{0.7}{$
\begin{array}{cccc|c}  
\rosalpha\SS & 0   & \cdots & 0 &\rosalpha\SS \\
0 & \sfrac{1}{\mratio} \rosalpha\FF                   & \cdots & 0 & \rosalpha\FSL[1]  \\
\vdots      &\vdots           & \ddots &  \vdots & \vdots  \\
0 & \sfrac{1}{\mratio}\one  \b\F* &   \ldots &  \sfrac{1}{\mratio} \rosalpha\FF &\rosalpha\FSL[\mratio] \\
\hline 
\rosalpha\SS & 0 &   \cdots &  0 & \rosalpha\SS   \\   
\Xhline{1.5pt}
0 &  \sfrac{1}{\mratio} \b\F* &   \ldots &  \sfrac{1}{\mratio} \b\F* & \b\S*  
\end{array}
$}
\raisebox{-15pt}{$\,,$}
\qquad
\raisebox{8.5pt}{
\scalebox{0.7}{$
\begin{array}{cccc|c}  
\rosgamma \SS & 0   & \cdots & 0 & \rosgamma \SS \\
0 & \sfrac{1}{\mratio} \rosgamma\FF                       & \cdots & 0 & \rosgamma\FSL[1]  \\
\vdots  & \vdots   & \ddots & \vdots  & \vdots  \\
0 & 0 &  \ldots &  \sfrac{1}{\mratio} \rosgamma\FF &\rosgamma\FSL[\mratio] \\
\hline 
\rosgamma \SS &  0 &   \cdots &  0 & \rosgamma \SS  
 \end{array}
 $}}
\raisebox{-15pt}{$\,.$}
\end{equation}
Step-predictor-corrector methods will be discussed in detail in \Cref{sec:Step-predictor-corrector}.

\section{Compound-first-step MR-GARK-ROS/ROW schemes}
\label{sec:compound-first-step}
%
Consider a telescopic compound-first-step method \eqref{eqn:mrGARK-ROW-butcher-compound-step} where the fast and the slow base methods coincide: $\rosalpha\FF=\rosalpha\SS=\rosalpha$, $\rosgamma\FF=\rosgamma\SS=\rosgamma$, $s\F=s\S=s$, and $\b\F=\b\S=\b$.

A natural choice for the fast/slow coupling coefficients is: 
\begin{subequations}
\label{eqn:Telescopic-ROS-coupling}
\begin{equation}
\label{eqn:Telescopic-ROS-FS}
\begin{split}
\rosgamma\FSL[\lambda] = \sfrac{1}{\mratio} \rosgamma + 
\sfrac{1}{\mratio} \widetilde{\rosgamma}\FSL[\lambda], \quad
\rosalpha\FSL[\lambda] = \sfrac{1}{\mratio}\,\rosalpha +
\sfrac{1}{\mratio}\,\widetilde{\rosalpha}\FSL[\lambda] +
\sfrac{1}{\mratio}\,\mathbf{F}(\lambda), 
\end{split}
\end{equation}
where we allow the additional coefficient matrices $\mathbf{F}(\lambda)$, $\widetilde{\rosalpha}\FSL[\lambda]$, and $\widetilde{\rosgamma}\FSL[\lambda]$ for more flexibility.  Here $\widetilde{\rosalpha}\FSL[\lambda]$ is assumed to be strictly lower triangular, and  $\widetilde{\rosgamma}\FSL[\lambda]$ lower triangular, which adds $\mratio\,(2s-3)$ degrees of freedom to coupling coefficients \eqref{eqn:Telescopic-ROS-FS}.
For internal consistency \eqref{eqn:internal-consistency-multirate} we ask that
\begin{equation}
\label{eqn:F-CFS-condition-IC}
\mathbf{F}(\lambda)\, \one = (\lambda-1)\,\one, \quad
\widetilde{\rosalpha}\FSL[\lambda]\,\one=\widetilde{\rosgamma}\FSL[\lambda]\,\one=0.
\end{equation}

A natural choice of  slow/fast coupling for a compound  first step is: 
\begin{equation}
\label{eqn:Telescopic-ROS-SF}
\rosalpha\SFL[1] = 
\mratio\,\rosalpha+\mratio\,\widetilde{\rosalpha}\SFL[1],  \qquad
\rosgamma\SFL[1] = 
\mratio\,\rosgamma+\mratio\,\widetilde{\rosgamma}\SFL[1],
\end{equation}
\end{subequations}
with $\widetilde{\rosalpha}\SFL[1]$ strictly lower triangular  and 
$\widetilde{\rosgamma}\SFL[1]$ lower triangular, where for internal consistency we impose
\begin{align*}
 \widetilde{\rosalpha}\SFL[1]\,\one  =
 \widetilde{\rosgamma}\SFL[1]\,\one = 0.
\end{align*}
This choice adds  $2s-3$ degrees of freedom to the coupling coefficients \eqref{eqn:Telescopic-ROS-SF}.

The particular choice of coupling coefficients \eqref{eqn:Telescopic-ROS-coupling} implies that:
\begin{equation*}
\begin{split}
&\rosalpha\FSL[\Sigma_0] = \rosalpha +
\sfrac{1}{\mratio}\,
\rosalpha\FSL[\Sigma_0]
+ 
\sfrac{1}{\mratio}\, \mathbf{F}(\Sigma_0),  \quad
\rosgamma\FSL[\Sigma_0] = \rosgamma + 
\sfrac{1}{\mratio}\,
\rosgamma\FSL[\Sigma_0], \\
&\rosalpha\SFL[\Sigma_0] = \rosalpha\SFL[1], \quad
\rosgamma\SFL[\Sigma_0] = \rosgamma\SFL[1],
\quad \rosalpha\SFL[\Sigma_1] = \rosgamma\SFL[\Sigma_1]  = 0,
\end{split}
\end{equation*}
where we have used the abbreviation
\begin{equation}
\label{eqn:F-sum-notation}
\mathbf{F}(\Sigma_k) \coloneqq \sum_{\lambda=1}^\mratio (\lambda-1)^k \,\mathbf{F}(\lambda).
\end{equation}
Assuming that the base ROS scheme has order three, the remaining MR-GARK-ROS order three conditions \eqref{eqn:MGARK-ROS-order3}  read:
\begin{eqnarray*}
	\rosb^T \, \rosbeta\SFL[1]\,\rose &=& \sfrac{\mratio^2}{6}, \\
	\rosb^T \,\left( \mathbf{F}(\Sigma_0) + 
\rosbeta\FSL[\Sigma_0] \right)\, \rose &=& \sfrac{\mratio(\mratio-1)}{6}.
\end{eqnarray*} 
Assuming that the base ROW scheme has order three, the remaining MR-GARK-ROW order three conditions \eqref{eqn:MGARK-ROW-order3} are:
\begin{align*}
\rosb^T \,  \rosalpha\SFL[1]\,\rosc  = \sfrac{\mratio^2}{6}, \quad
\rosb^T\,\left(\mathbf{F}(\Sigma_0) + \boldsymbol{\tilde \alpha}\FSL[\Sigma_0] \right)\, \rosc = \sfrac{\mratio(\mratio-1)}{6}, \quad
\rosb^T \, \rosgamma\SFL[1]\,\rosc = 0, \\
\rosb^T \, \left( \rosgamma + \sfrac{1}{\mratio} \widetilde{\rosbeta}\FSL[\Sigma_0] 
\right) \,\rosc = 0, \quad
\rosb^T \,  \rosalpha\SFL[1] \, \rosg = 0, \quad
b^T\,\left( \mathbf{F}(\Sigma_0) +  \boldsymbol{\tilde \alpha}\FSL[\Sigma_0]  \right) \, \rosg = 0, \\
\rosb^T \,  \rosgamma\SFL[1] \, \rosg = 0, \quad
\rosb^T \, \left(  \rosgamma + \sfrac{1}{\mratio} \widetilde{\rosbeta}\FSL[\Sigma_0] \right) \, \rosg = 0.
\end{align*}
%
%
%
Assuming that the base ROS scheme has order three, the remaining MR-GARK-ROS order four conditions \eqref{eqn:MGARK-ROS-order4} are:
\begin{eqnarray*}
\rosb^T \, \Big(  \rosalpha\SFL[1]\,\rose   \times \rosc \Big) &=& \sfrac{\mratio^2}{8},  \\
%
\rosb^T \, \Big(  (\mathbf{F}(\Sigma_1)+\boldsymbol{\tilde \alpha}\FSL[\Sigma_1])\,\rose + \big( (\mathbf{F}(\Sigma_0)+\boldsymbol{\tilde \alpha}\FSL[\Sigma_0])\, \rose\big) \times \rosc \Big) &=& 
 \sfrac{\mratio(\mratio-1)}{8}\, \left(\mratio+1 -  4\,\rosb^T \rosalpha\,\rose\right), \\
%
\rosb^T \,  \rosbeta\SFL[1]\,\rosc^{\times 2} &=& \sfrac{\mratio^3}{12}, \\
%
\rosb^T \, (\mathbf{F}(\Sigma_0)+\widetilde{\rosbeta}\FSL[\Sigma_0]) \, \rosc^{\times 2} &=& \sfrac{\mratio(\mratio-1)}{12}, \\
%
\rosb^T \,\rosbeta\, \rosbeta\SFL[1]\,\rose &=&  \sfrac{\mratio^2}{24},  \\
%
\rosb^T \,\rosbeta\SFL[1]\,\rosbeta \,\rose  &=&  \sfrac{\mratio^3}{24},  \\
%
\rosb^T \,\rosbeta\SFL[1]\, (\mathbf{F}(1) + \widetilde{\rosbeta}\FSL[1])\,\rose &=& \sfrac{\mratio^2 (1-\mratio)}{24}, \\
%
\rosb^T  \,\rosbeta\,(\mathbf{F}(\Sigma_0)+\widetilde{\rosbeta}\FSL[\Sigma_0])\,\rose &=&  \sfrac{\mratio(\mratio^2-1)}{24},  \\
%
\rosb^T \, (\mathbf{F}(\Sigma_0)+\widetilde{\rosbeta}\FSL[\Sigma_0])\,\rosbeta\,\rose &=& \sfrac{\mratio(\mratio-1)}{24},  \\
%
\rosb^T \, \big(\mathbf{F}(1)+\widetilde{\rosbeta}\FSL[1] \big) \,\rosbeta\SFL[1]\, \rose &=& \sfrac{\mratio^2(\mratio-1)}{24}.
\end{eqnarray*}

\begin{remark}
Order four coupling conditions for compound-first-step MR-GARK-ROW schemes can be derived in an analogous manner from the general order conditions \eqref{eqn:GARK-ROW-order4-conditions}.
\end{remark}

Using the framework of compound-first-step schemes we derive embedded  MR-GARK-ROS methods of order (2)3 (main method of order three, with embedded scheme of order two)  which fulfill the time-lagged Jacobian order conditions.
\begin{example}[Implicit-implicit case]
\label{ex.imim}
An embedded implicit-implicit MR-GARK-ROS scheme of order (2)3,
which automatically fulfills the time-lagged Jacobian order conditions~\eqref{eqn:GARK-ROW-time-lagged-order3-conditions} due~\Cref{remark.time-lagged}, with $\gamma$ and $\beta_{2,1}$ as free parameters,
is given by 
\begin{align*}
 \rosb & = 
    \begin{bmatrix}
    \frac{1}{6} & \frac{4}{6} & \frac{1}{6}
    \end{bmatrix}^T, \qquad  
    \hat{\rosb} = \begin{bmatrix}
  {\scriptstyle 1}- \frac{1-2\,\gamma}{2\,\beta_{2,1}}  & \frac{1-2\,\gamma}{2\,\beta_{2,1}} & {\scriptstyle 0}
    \end{bmatrix}^T, \\
    \rosalpha & =
    \begin{bmatrix}
    \hphantom{{\scriptstyle -}}  {\scriptstyle 0} & {\scriptstyle 0} & {\scriptstyle 0} \\  \hphantom{{\scriptstyle -}} \frac{1}{2} & {\scriptstyle 0} & {\scriptstyle 0}  \\{\scriptstyle -1} & {\scriptstyle 2} & {\scriptstyle 0}
    \end{bmatrix}, \quad
    \rosbeta = 
  \begin{bmatrix}
    {\scriptstyle \gamma} & {\scriptstyle 0} & {\scriptstyle 0} \\ 
    {\scriptstyle \beta_{2,1}} & {\scriptstyle \gamma} & {\scriptstyle 0}  \\ 
    {\scriptstyle 3-6 \gamma-4 \beta_{2,1} -} \frac{6(\gamma^2+\gamma)-1}{\beta_{2,1}} & \frac{6(\gamma^2-\gamma)+1}{\beta_{2,1}}  & {\scriptstyle \gamma}
    \end{bmatrix},  \\
    \rosalpha\FSL[\lambda] & = 
    \sfrac{1}{\mratio} \Bigl( \rosalpha+ \mathbf{F}(\lambda) \Bigr), \quad  \rosbeta\FSL[\lambda]=
   \sfrac{1}{\mratio} \left( \rosbeta + \hat{\rosbeta}  + \mathbf{F}(\lambda) \right), \\  
      \rosalpha\SFL[1]  & = \mratio\, \rosalpha, \quad
    \rosbeta\SFL[1]  = 
    \mratio\, \left( \rosbeta + \hat{\rosbeta} \right), \\ 
 \hat{\rosbeta}  & = \begin{bmatrix}
    {\scriptstyle 0} & {\scriptstyle 0} & {\scriptstyle 0} \\ 
    {\scriptstyle 0} & {\scriptstyle 0} & {\scriptstyle 0} \\
    {\scriptstyle -\hat{\beta}}  & {\scriptstyle \hat{\beta}}  & {\scriptstyle 0}
        \end{bmatrix}, \quad
  \hat{\beta}  = \sfrac{\mratio-1}{\beta_{2,1}}, \\
        \mathbf{F}(\lambda) & = (\lambda-1)
        \begin{bmatrix}
{\scriptstyle \hat \alpha}  & {\scriptstyle 1- \hat \alpha}  & {\scriptstyle 0} \\   
{\scriptstyle \hat \alpha}  & {\scriptstyle 1- \hat \alpha}  & {\scriptstyle 0} \\ 
{\scriptstyle \hat \alpha}  & {\scriptstyle 1- \hat \alpha}  & {\scriptstyle 0} \\ 
        \end{bmatrix}, 
        \quad \hat \alpha = \sfrac{\beta_{2,1}+\gamma}{\beta_{2,1}}.
\end{align*}
\end{example}

\begin{example}[Implicit-explicit case]
For the implicit-explicit case we select $\rosbeta\FFL[\lambda]=\rosalpha\FF=\rosalpha$ for 
$\lambda=2,\ldots,\mratio$, i.e., the base scheme for the last $\mratio-1$ steps of the fast part is explicit. To obtain a fully-implicit compound step we set $\rosbeta\FFL[1]=\rosbeta$. We choose $\rosalpha, \rosbeta, \rosb, \hat \rosb$ and the other coupling coefficients as in \Cref{ex.imim} above, but replace $\hat \alpha$ with:
\begin{align*}
    \hat \alpha =
    \frac{3(\mratio+1)(\beta_{2,1}+\gamma)-\beta_{2,1}-\mratio}{3 \mratio\, \beta_{2,1}}.
\end{align*}.
\end{example}

\section{Step-predictor-corrector methods}
\label{sec:Step-predictor-corrector}

%
%
We consider step predictor-corrector (SPC) methods described in section \ref{subsec:step-predictor-corrector}. The computations associated with the Butcher tableau  \eqref{eqn:MROS/MROW-SPC-butcher} proceed as follows. 
First, the ``predictor'' applies the slow base scheme $(\rosb\S,\rosalpha\SS,\rosgamma\SS)$ over a macro-step $H$ to solve the entire coupled system:
\begin{subequations}
\label{eqn:MROS/MROW-SPC}
\begin{eqnarray}
\mathbf{k} &=& H\,(\fun\F+\fun\S)\left( \one_{s\S} \otimes \y_{n-1}  + 
 \rosalpha\SS \kron{\nvar} \mathbf{k}\right)  \nonumber \\
\label{eqn:MROS/MROW-Kfull-CS}
& &  +  \left( \Idstage \otimes H (\Lb\F+ \Lb\S) \right) \,
 \left( \rosgamma\SS  \kron{\nvar} \mathbf{k} \right),
\end{eqnarray}
which gives the following values of the slow stages:
\begin{eqnarray}
\label{eqn:MROS/MROW-KS-CS}
\quad \mathbf{k}\S &=& H\,\fun\S\left( \one_{s\S} \otimes \y_{n-1}  + 
 \rosalpha\SS \kron{\nvar} \mathbf{k} \right)    +  \bigl( \Idstage \otimes H \Lb\S) \,\left( \rosgamma\SS  \kron{\nvar} \mathbf{k}  \right).
\end{eqnarray}
\end{subequations} 
Next, the ``corrector'' applies the fast base scheme $(\rosb\F,\rosalpha\FF,\rosgamma\FF)$ over $\mratio$ micro-steps of size $h$. The fast stages $\mathbf{k}\FL[\lambda]$ are computed using formula \eqref{eqn:MROS/MROW-KF}, and the next step solution $\y_{n}$ using formula \eqref{eqn:MROS/MROW-Sol}. Note that the coupling matrices $\rosalpha\FSL[\lambda]$ and $\rosgamma\FSL[\lambda]$ do not need to be triangular, and can have any fill-in structure, since all $\mathbf{k}\S$ are known before the start of the micro-steps.

\subsection{Order conditions}

\paragraph{Internal consistency}
\ifnum \value{book}=1
The vectors \eqref{eqn:terminal-vectors-multirate} are:
\begin{equation}
\begin{split}
\mathbf{c}\SS & = \mathbf{c}\SF = \rosc\SS, \\
\mathbf{c}\FS & = \Big[\rosc\SS; \big[ \rosc\FSL[\lambda] \big]_{1 \le \lambda \le \mratio} \Big], \\
\mathbf{c}\FF & = \Big[\rosc\SS; \big[ \sfrac{\lambda-1}{\mratio}\, \one\F +  \sfrac{1}{\mratio}\, \rosc\FF \big]_{1 \le \lambda \le \mratio}\Big], \\
\mathbf{g}\SS & = \mathbf{g}\SF = \rosg\SS, \\
\mathbf{g}\FS & = \Big[\rosg\SS; \big[ \rosg\FSL[\lambda] \big]_{1 \le \lambda \le \mratio} \Big], \\
\mathbf{g}\FF & = \Big[\rosg\SS; \big[ \sfrac{1}{\mratio}\, \rosg\FF \big]_{1 \le \lambda \le \mratio} \Big].
\end{split}
\end{equation}
\fi
The internal consistency conditions \eqref{eqn:ROW-internal-consistency} read:
\begin{equation}
\label{eqn:SPC-internal-consistency}
\begin{split}
& 
\rosc\FSL[\lambda] = \sfrac{\lambda-1}{\mratio}\, \one\F +  \sfrac{1}{\mratio}\, \rosc\FF, \quad \lambda = 1, \dots,  \mratio, \\
&
\rosg\FSL[\lambda] = \sfrac{1}{\mratio}\, \rosg\FF, \quad \lambda = 1, \dots,  \mratio.
\end{split}
\end{equation}

The accuracy analysis in this section assumes internally consistent SPC-MR-GARK-ROS/ROW methods \eqref{eqn:SPC-internal-consistency}, where the slow and fast base methods are ROS/ROW schemes of the corresponding order. Therefore the analysis below focuses only on the remaining coupling conditions. 

\paragraph{Second order conditions}
\ifnum \value{book}=1
The second order ROS coupling condition \eqref{eqn:GARK-ROS-order2-conditions} is:
\begin{align}
\label{eqn:SPC-ROS-order2}
\rosb\F* \, \rose\FSL[\Sigma_0]  = \sfrac{\mratio}{2},
\end{align}
and is satisfied if the method is internally consistent \eqref{eqn:SPC-internal-consistency}.

The second order ROW coupling order conditions \eqref{eqn:GARK-ROW-order2-conditions} are:
\begin{align}
\label{eqn:SPC-ROW-order2}
\rosb\F* \, \rosc\FSL[\Sigma_0]  = \sfrac{\mratio}{2}, \qquad \rosb\F* \, \rosg\FSL[\Sigma_0]  = 0,
\end{align}
and are both satisfied if the method is internally consistent \eqref{eqn:SPC-internal-consistency}.
\else
For internally consistent SPC-MR-GARK-ROS/ROW methods the coupling conditions \eqref{eqn:GARK-ROS-order2-conditions} and \eqref{eqn:GARK-ROW-order2-conditions}, respectively, are automatically satisfied.
\fi

\paragraph{Third order conditions}
For SPC-MR-GARK-ROS methods there is a single remaining order three coupling condition \eqref{eqn:GARK-ROS-order3-conditions}, as follows:
\begin{align}
\label{eqn:SPC-ROS-order3}
\rosb\F* \, \rosbeta\FSL[\Sigma_0] \, \rose\SS = \sfrac{\mratio}{6}. 
\end{align}
The third order conditions for time-lagged Jacobians \eqref{eqn:GARK-ROW-time-lagged-order3-conditions} are automatically satisfied due to internal consistency.

For SPC-MR-GARK-ROW methods the order three coupling conditions \eqref{eqn:GARK-ROW-order3-conditions} are:
\begin{align}
\label{eqn:SPC-ROW-order3}
\rosb\F* \, \rosalpha\FSL[\Sigma_0]\, \rosc\SS &= \sfrac{\mratio}{6}, &
\rosb\F* \, \rosgamma\FSL[\Sigma_0]\,\rosc\SS &= 0, \\
\nonumber
\rosb\F* \, \rosalpha\FSL[\Sigma_0]\,\rosg\SS &= 0, &
\rosb\F* \, \rosgamma\FSL[\Sigma_0]\, \rosg\SS &= 0.
\end{align}

\paragraph{Fourth order conditions}
For SPC-MR-GARK-ROS methods \eqref{eqn:SPC-internal-consistency} there are four remaining order four coupling ROS conditions \eqref{eqn:GARK-ROS-order4-conditions}, as follows:
\begin{subequations}
\label{eqn:SPC-ROS-order4}
\begin{align}
\label{eqn:SPC-ROS-order4-a}
\rosb\F* \,\left( \rosalpha\FSL[\Sigma_1] \,\rose\SS  + (\rosalpha\FSL[\Sigma_0] \,\rose\SS) \times \rosc\FF \right)  &= \sfrac{\mratio^2}{8}, \\%
%
\label{eqn:SPC-ROS-order4-b}
\rosb\F*\,  \rosbeta\FSL[\Sigma_0] \, \rosc\SS\,^{\times 2}   &= \sfrac{\mratio}{12}, \\ 
%
%
\label{eqn:SPC-ROS-order4-c}
\rosb\F*\,   \left( \sum_{\lambda=1}^\mratio\sum_{i=1}^{\lambda-1} \rosbeta\FSL[i]  + \rosbeta\FF\,\rosbeta\FSL[\Sigma_0] \right)\,\rose\SS 
&=  \sfrac{\mratio^2}{24}, \\
%
%
\label{eqn:SPC-ROS-order4-d}
\rosb\F*\, \rosbeta\FSL[\Sigma_0] \,\rosbeta\SS\,\rose\SS &=  \sfrac{\mratio}{24}.
%
%
%
\end{align}
\end{subequations}

\begin{remark}
Order four coupling conditions for SPC-MR-GARK-ROW schemes can be derived in an analogous manner from the general order conditions \eqref{eqn:GARK-ROW-order4-conditions}.
\end{remark}

%

\subsection{Telescopic SPC methods}
%
Consider a telescopic SPC method \eqref{eqn:MROS/MROW-SPC-butcher} where the fast and the slow base methods coincide, and are both equal to  $(\rosb,\rosalpha,\rosgamma)$; the vectors \eqref{eqn:terminal-vectors} are $\rosc = \rosalpha\,\one$, $\rosg = \rosgamma\,\one$, and $\rose = \rosc + \rosg$.

The internal consistency equations \eqref{eqn:SPC-internal-consistency} read:
\begin{equation*}
\rosc\FSL[\lambda] = \sfrac{\lambda-1}{\mratio}\, \one +  \sfrac{1}{\mratio}\, \rosc, \quad 
\rose\FSL[\lambda] = \sfrac{\lambda-1}{\mratio}\, \one +  \sfrac{1}{\mratio}\, \rose, \quad
\rosg\FSL[\lambda] = \sfrac{1}{\mratio}\, g,
\quad \lambda = 1, \dots,  \mratio.
\end{equation*}
We consider the following natural choice for the coupling coefficients:
\begin{subequations}
\label{eqn:SPC-ROS-coupling}
\begin{equation}
\label{eqn:F-coupling}
\rosgamma\FSL[\lambda] = \sfrac{1}{\mratio} \rosgamma, \quad
\rosalpha\FSL[\lambda] = \sfrac{1}{\mratio}\,\rosalpha + \sfrac{1}{\mratio}\,\mathbf{F}(\lambda), \quad
\rosbeta\FSL[\lambda] = \sfrac{1}{\mratio}\,\rosbeta + \sfrac{1}{\mratio}\,\mathbf{F}(\lambda),
\end{equation}
where, in order to ensure internal consistency, we ask that
\begin{equation}
\label{eqn:F-condition-IC}
\mathbf{F}(\lambda)\, \one = (\lambda-1)\,\one, \quad \lambda = 1,\dots, \mratio.
\end{equation}
\end{subequations}
The second order coupling conditions for SPC-MR-GARK-ROS/ROW schemes are automatically satisfied due to internal consistency.

Using notation \eqref{eqn:F-sum-notation} and the coupling coefficients \eqref{eqn:F-coupling} we have:
\begin{equation}
\label{eqn:F-coupling-sigmas}
\begin{split}
\rosgamma\FSL[\Sigma_0] &= \rosgamma, \quad
\rosgamma\FSL[\Sigma_1] = \sfrac{\mratio-1}{2} \rosgamma, \\
\rosalpha\FSL[\Sigma_0] &= \rosalpha + \sfrac{1}{\mratio}\,\mathbf{F}(\Sigma_0), \quad
\rosalpha\FSL[\Sigma_1] = \sfrac{\mratio-1}{2} \rosalpha + \sfrac{1}{\mratio}\,\mathbf{F}(\Sigma_1).
\end{split}
\end{equation}
\paragraph{Third order conditions}
The third order SPC-MR-GARK-ROS coupling condition \eqref{eqn:SPC-ROS-order3} reads:
\begin{equation}
\label{eqn:SPCTEL-ROS-order3}
\rosb^T \, \mathbf{F}(\Sigma_0)\, \rose = \sfrac{\mratio\,(\mratio-1)}{6}. 
\end{equation}
The third order SPC-MR-GARK-ROW coupling conditions \eqref{eqn:SPC-ROW-order3} are:
\begin{equation}
\label{eqn:SPCTEL-ROW-order3}
\rosb^T \,  \mathbf{F}(\Sigma_0) \, \rosc = \sfrac{\mratio\,(\mratio-1)}{6}, \qquad
\rosb^T \, \mathbf{F}(\Sigma_0)  \,\rosg= 0.
\end{equation}
%
%
The following choice of rank-one coupling matrix ensures third order: 
\begin{align*}
\mathbf{F}(\lambda) &=  (\lambda-1)\, \one\,\mathbf{v}_1^T, 
\end{align*}
where impose the internal consistency \eqref{eqn:F-condition-IC}, the coupling condition \eqref{eqn:SPCTEL-ROS-order3}, and the second condition \eqref{eqn:SPCTEL-ROW-order3} by the following equations, respectively:
\begin{align*}
\mathbf{v}_1^T\, \one = 1, \quad \mathbf{v}_1^T\, \rose = \sfrac{1}{3}, \quad \mathbf{v}_1^T\, \rosg = 0.
\end{align*}
If the base method is a third order ROS scheme, then $\mathbf{v}_1^T = 2\,\rosb^T\, \rosbeta$ offers a solution of these equations.

\paragraph{Fourth order conditions}
The SPC-MR-GARK-ROS order four coupling conditions \eqref{eqn:SPC-ROS-order4-b} and \eqref{eqn:SPC-ROS-order4-d} are automatically satisfied by the choice of coupling coefficients \eqref{eqn:SPC-ROS-coupling}. The remaining order four coupling conditions \eqref{eqn:SPC-ROS-order4-a} and \eqref{eqn:SPC-ROS-order4-c} are:
\begin{equation}
\label{eqn:SPCTEL-ROS-order4}
\begin{split}
\rosb^T \,\left( \mathbf{f}(\Sigma_1)  + \mathbf{f}(\Sigma_0) \times  c \right)  &= \sfrac{\mratio\,(\mratio-1)(\mratio + 1 - 4\,\rosb^T \,\rosalpha \,\rose)}{8}, \\
\rosb^T \,  \left(  \sum_{\lambda=1}^\mratio\sum_{i=1}^{\lambda-1} \mathbf{f}(i)  + \rosbeta\,\mathbf{f}(\Sigma_0)  \right)
&=  \sfrac{\mratio\,(\mratio-1)^2}{24}, \\
\textnormal{where} \qquad \mathbf{f}(\lambda) \coloneqq \mathbf{F}(\lambda)\,\rose, \quad
\mathbf{f}(\Sigma_0) &\coloneqq \mathbf{F}(\Sigma_0)\,\rose, \quad
\mathbf{f}(\Sigma_1) \coloneqq \mathbf{F}(\Sigma_1)\,\rose.
\end{split}
\end{equation}
\paragraph{Polynomial coupling matrix}
We build the coupling matrix $\mathbf{F}(\lambda)$ as a polynomial with matrix coefficients:
\begin{align*}
\mathbf{F}(\lambda) &=  \sum_{i = 1}^K (\lambda-1)^i\, \mathbf{D}_i.
\end{align*}
Let $K=2$. For internal consistency \eqref{eqn:F-condition-IC} we require
\[
 \mathbf{D}_1\,\one = \one, \qquad \mathbf{D}_2\,\one = \zero.
\]
The third order ROS coupling conditions \eqref{eqn:SPCTEL-ROS-order3} are: 
\begin{align*}
\rosb^T\, \mathbf{D}_1\, \rose = \frac{1}{3}, \qquad \rosb^T\, \mathbf{D}_2\, \rose = 0,
\end{align*}
and the third order ROW coupling conditions \eqref{eqn:SPCTEL-ROW-order3} read:
\begin{align*}
\rosb^T\, \mathbf{D}_1\, \rosc = \frac{1}{3}, \qquad \rosb^T\, \mathbf{D}_2\, \rosc = 0, \qquad
\rosb^T\, \mathbf{D}_1\, \rosg = 0, \qquad \rosb^T\, \mathbf{D}_2\, \rosg = 0.
\end{align*}
The fourth order ROS coupling conditions \eqref{eqn:SPCTEL-ROS-order4} are:
\begin{align*}
\rosb^T\, (\mathbf{D}_1\, \rose \times \rosc) &= \frac{3}{8} - \rosb^T \rosalpha\, \rose, & \rosb^T\, (\mathbf{D}_2\, \rose \times \rosc) &= \sfrac{1}{24}, \\
\rosb^T\, \rosbeta\, \mathbf{D}_1\, \rose &= \sfrac{1}{8}, &  \rosb^T\, \rosbeta\, \mathbf{D}_2\, \rose &= -\sfrac{1}{24}.
\end{align*}

\section{Multirate Infinitesimal Step Methods}
\label{sec:infinitesimal-step}

We now consider MR-GARK-ROS/ROW methods where the micro-steps can be arbitrarily small. We call these methods multirate ``infinitesimal step'', or 
MRI-GARK-ROS/ROW for short; they offer extreme flexibility since they allow to solve the fast sybsystem with any sufficiently accurate discretization method and sequence of step sizes.

\begin{definition}[MRI-GARK-ROS/ROW methods]
Consider a base slow ROS/ROW scheme $(\rosb\S,\rosalpha\SS,\rosgamma\SS)$ with non-decreasing abscissae $\rosc\SS_1 \le \rosc\SS_2 \le \cdots \le \rosc\SS_{s\S}$, and denote:
\begin{eqnarray*}
\rosc\SS_0=0; \qquad
\Delta \rosc\SS_i \coloneqq \rosc\SS_{i}-\rosc\SS_{i-1}, \quad i=1,\dots,s.
\end{eqnarray*}
A multirate infinitesimal step GARK ROS/ROW method  advances the solution of the fast-slow partitioned system \eqref{eqn:multirate-ode} via the following computational process:
\begin{subequations}
\label{eqn:MROS/MROW-INF}
\begin{eqnarray}
&&\widetilde{\y}_{n-1}  = \y_{n-1}; \\[3pt]
%
%
&&\textnormal{For }\lambda=1,\dots,s\S: \nonumber \\
\label{eqn:MROS/MROW-INF-fast}
&&
\begin{cases}
\mathbf{v}_\lambda(0) = \widetilde{\y}_{n-1+(\lambda-1)/s}, \\
\mathbf{v}_\lambda'(\theta) = \Delta \rosc\SS_\lambda\,\fun\F\Bigl(t_{n-1}+\rosc\SS_{\lambda-1}\,H + \Delta \rosc\SS_\lambda\theta, \\
\qquad\qquad\qquad \mathbf{v}_\lambda(\theta)  +  \sum_{j=1}^{\lambda-1} \mathbf{r}_{\lambda,j} ({\scriptstyle \frac{\theta}{H}})  \, \mathbf{k}\S_j \Bigr), 
\qquad \theta \in [0,H],  \\[3pt]
\widetilde{\y}_{n-1+\lambda/s}  = \mathbf{v}_\lambda(H); 
 \end{cases}
 \\[3pt]
%
\label{eqn:MROS/MROW-INF-slow}
&&\mathbf{k}\S_\lambda = H\,\fun\S\left( \widetilde{\y}_{n-1+\lambda/\mratio} 
+   \sum_{j=1}^{\lambda-1} \rosalpha\SS_{\lambda,j} \mathbf{k}\S_j \right)  \\
&&\qquad   +  H\, \Lb\S \,\left(
\sum_{\ell=1}^{\lambda} \mathbf{q}_{\lambda,\ell}\,  (\widetilde{\y}_{n-1+\ell/\mratio}-\widetilde{\y}_{n-1+(\ell-1)/\mratio})
+ \sum_{j=1}^{\lambda}\rosgamma\SS_{\lambda,j} \mathbf{k}\S_j  \right); \nonumber  \\
\label{eqn:MROS/MROW-INF-solution}
&& \y_{n} = \widetilde{\y}_{n}  +  \rosb\S*  \kron{\nvar} \mathbf{k}\S.
\end{eqnarray}
\end{subequations} 
A modified fast ODE \eqref{eqn:MROS/MROW-INF-fast} is integrated between consecutive stages of the base slow method. The slow components influence the fast dynamics via the time dependent coefficients $\mathbf{r}_{\lambda}(\cdot) \in \Re^{s\S}$, $\lambda = 1,\dots,s\S$. The fast solutions impact the computation of the slow stages \eqref{eqn:MROS/MROW-INF-slow} via the coupling coefficients $\mathbf{q}_{\lambda} \in \Re^{s\S}$, $\lambda = 1,\dots,s\S$. The next step solution \eqref{eqn:MROS/MROW-INF-solution} combines the fast solution $\widetilde{\y}_{n}$ and the slow solution increment given by the stages $\mathbf{k}\S$.
\end{definition}

To analyze the scheme \eqref{eqn:MROS/MROW-INF} we start by discretizing each modified fast ODE \eqref{eqn:MROS/MROW-INF-fast} with an explicit Runge-Kutta scheme $(\rosb\F,\rosalpha\FF)$ of arbitrary accuracy:
\begin{equation}
\label{eqn:MROS/MROW-INF-fast-RK}
\begin{split}
& \mathbf{k}_i\FL[\lambda] = H \,\fun\F\Bigl( \widetilde{\y}_{n-1+(\lambda-1)/\mratio}  + 
\sum_{j=1}^{i-1} \rosalpha\FF_{i,j}  \,\Delta \rosc\SS_\lambda \mathbf{k}_j\FL[\lambda]  + 
 \sum_{j=1}^{\lambda-1} \rosalpha\FSL[\lambda]_{i,j}  \, \mathbf{k}\S_j \Bigr),\\
 & \qquad\qquad  i = 1,\dots,s\F, \\
& \widetilde{\y}_{n-1+\lambda/\mratio} = \widetilde{\y}_{n-1+(\lambda-1)/\mratio}  + 
 \sum_{i=1}^{s\F} \rosb\F_i  \,\Delta \rosc\SS_\lambda \mathbf{k}\FL[\lambda]_i. 
 \end{split}
\end{equation}
The fast discrete stages \eqref{eqn:MROS/MROW-INF-fast-RK}, together with the slow stages \eqref{eqn:MROS/MROW-INF-slow} and the next step solution \eqref{eqn:MROS/MROW-INF-solution}, form an IMEX GARK ROS/ROW method \eqref{eqn:MROS/MROW} with $\mratio = s\S$, as described in Section \ref{subsec:coupling-imex}. The slow/fast coupling coefficients have the following particular structure:
\begin{eqnarray*}
\rosalpha\FSL[\lambda]_{i,j} &\coloneqq& \mathbf{r}_{\lambda,j} (\rosc\F_i), \\
\rosalpha\SFL[\lambda] &\coloneqq& \mathbf{p}_\lambda\,\rosb\F*, \quad
 \mathbf{p}_\lambda \coloneqq \begin{bmatrix} \zero_{(\lambda-1) \times 1} \\ \one_{(s\S+1-\lambda) \times 1} \end{bmatrix}, \\
 %
 \rosgamma\SFL[\lambda] &\coloneqq&  \boldsymbol{q}_\lambda\,\rosb\F*,
 \quad \mathbf{q}_\lambda \coloneqq \begin{bmatrix} \zero_{(\lambda-1) \times 1} \\ q_{\lambda} \end{bmatrix},
 \quad  q_{\lambda} \in \Re^{s\S+1-\lambda}.
\end{eqnarray*}

The Butcher tableau \eqref{eqn:mrGARK-ROW-butcher} of the resulting IMEX GARK ROS/ROW scheme is:
\begin{equation*}
\renewcommand{\arraystretch}{1.5}
\scalebox{0.7}{$
\begin{array}{ccc|cccc}  
 \Delta \rosc\SS_1 \rosalpha\FF                       & \cdots & 0 & \rosalpha\FSL[1]  \\
\vdots                      & \ddots &   & \vdots  \\
 \Delta \rosc\SS_1 \one  \rosb\F*   & \ldots &  \Delta \rosc\SS_{s\S}  \rosalpha\FF &\rosalpha\FSL[s\S] \\
\hline 
 \Delta \rosc\SS_1 \mathbf{p}_1 \rosb\F*  & \cdots &  \Delta \rosc\SS_{s\S} \mathbf{p}_{s\S} \rosb\F*& \rosalpha\SS   \\   \Xhline{1.5pt} 
  \Delta \rosc\SS_1 \rosb\F*  & \ldots &  \Delta \rosc\SS_{s\S} \rosb\F* & \rosb\S*  
\end{array}
$},
\quad
\raisebox{6pt}{
\scalebox{0.7}{$
\begin{array}{ccc|cccc}  
 0                       & \cdots & 0 & 0 \\
\vdots                   & \ddots &   & \vdots  \\
 0   & \ldots &  0 & 0 \\
\hline 
 \Delta \rosc\SS_1 \mathbf{q}_1 \rosb\F* &    \cdots &  \Delta \rosc\SS_{s\S} \mathbf{q}_{s\S} \rosb\F*& \rosgamma\SS  
 \end{array}
 $}}.
\end{equation*}

\paragraph{Internal consistency}
The internal consistency conditions \eqref{eqn:internal-consistency-multirate-csf} and \eqref{eqn:internal-consistency-multirate-gfs} are automatically satisfied. Conditions \eqref{eqn:internal-consistency-multirate-cfs} and \eqref{eqn:internal-consistency-multirate-gsf} read:
\begin{equation*}
\begin{split}
& \rosc\FSL[\lambda] = \rosc\SS_{\lambda-1}\, \one\F +  \Delta \rosc\SS_\lambda\, \rosc\FF, \quad \lambda = 1, \dots, s, \\
& \sum_{\lambda=1}^{s\S} \Delta \rosc\SS_\lambda q_\lambda = \rosg\SS.
\end{split}
\end{equation*}

The following choice obeys the internal consistency condition:
\begin{eqnarray*}
\begin{split}
& \mathbf{r}_{\lambda}(\theta) = \sum_{k \ge 0} \mathbf{r}_{\lambda,k}\,\theta^k \in \Re^{s\S}, \quad
\rosalpha\FSL[\lambda] \coloneqq \sum_{k \ge 0} \rosc\F[\times k]\, \mathbf{r}_{\lambda,k}^T, \\
& \mathbf{r}_{\lambda,0}^T\, \one\S = \rosc\SS_{\lambda-1},
\quad 
\mathbf{r}_{\lambda,1}^T\, \one\S = \Delta \rosc\SS_\lambda, \quad
\mathbf{r}_{\lambda,k}^T\, \one\S = 0,~k \ge 2.
\end{split}
\end{eqnarray*}
We have that:
\begin{eqnarray*}
\rosb\F*\,\rosalpha\FSL[\lambda] = \bar{\mathbf{r}}_{\lambda}^T, \quad 
\bar{\mathbf{r}}_{\lambda} \coloneqq \sum_{k \ge 0} \sfrac{1}{k+1}\, \mathbf{r}_{\lambda,k}, \quad 
\hat{\mathbf{r}}_{\lambda} \coloneqq \sum_{k \ge 0} \sfrac{1}{k+2}\, \mathbf{r}_{\lambda,k}.
\end{eqnarray*}
\paragraph{Third order conditions}
The MRI-GARK-ROS coupling conditions \eqref{eqn:MGARK-ROS-order3} read: 
\begin{eqnarray*}
\rosb\S* \,\sum_{\lambda=1}^{s\S} \mathbf{q}_\lambda\, \Big(   \rosc\SS[2]_{\lambda} - \rosc\SS[2]_{\lambda-1} \Big) = 0,\qquad
%
\sum_{\lambda=1}^{s\S}   \Delta \rosc\SS_\lambda\, \bar{\mathbf{r}}_{\lambda}^T\, \rose\SS = \sfrac{1}{6}.
\end{eqnarray*} 
The order three MRI-GARK-ROW coupling conditions \eqref{eqn:MGARK-ROW-order3} are:
\begin{eqnarray*}
%
%
\rosb\S* \,  \sum_{\lambda=1}^s  \mathbf{q}_\lambda\, \Big( \rosc\SS[2]_{\lambda} - \rosc\SS[2]_{\lambda-1} \Big) &=& 0, \\
%
\sum_{\lambda=1}^s  \Delta \rosc\SS_\lambda\,\bar{\mathbf{r}}_{\lambda}^T\, \rosc\SS = \sfrac{1}{6}, \qquad
%
\sum_{\lambda=1}^s  \Delta \rosc\SS_\lambda\,\bar{\mathbf{r}}_{\lambda}^T\,\rosg\SS &=& 0.
\end{eqnarray*}

\ifnum \value{book}=1
\paragraph{Fourth order conditions}
The MR-GARK-ROS coupling conditions \eqref{eqn:MGARK-ROS-order3} are: 
\begin{eqnarray*}
\sum_{\lambda=1}^s \Delta \rosc\SS_\lambda\,  \Big( \rosc\SS_{\lambda-1}\,\bar{\mathbf{r}}_{\lambda}  +  \Delta \rosc\SS_\lambda\, \hat{\mathbf{r}}_{\lambda} \Big) ^T\,\rose\SS= \frac{1}{8}, \\
\sum_{\lambda=1}^s \big(\rosc\SS[2]_{\lambda} -  \rosc\SS[2]_{\lambda-1}\big)\, \sum_{i=\lambda}^s \rosb\S_i \,  \rosc\SS_i = \frac{1}{4}, \\
%
\sum_{\lambda=1}^s \Delta \rosc\SS_\lambda\, \bar{\mathbf{r}}_{\lambda}^T\, \rosc\SS[\times 2]   = \sfrac{1}{12}, \\ 
\rosb\S* \,\sum_{\lambda=1}^s  (\boldsymbol{p}_\lambda + \boldsymbol{q}_\lambda)\, (\rosc\SS[3]_{\lambda} -  
\rosc\SS[3]_{\lambda-1} ) = \sfrac{1}{4}, \\ 
%
\sum_{\lambda=1}^{s\S} \Delta \rosc\SS[2]_\lambda \, \sum_{k \ge 0} \sfrac{1 + (s-\lambda)(k + 3)}{(k + 1)(k + 3)} \, \mathbf{r}_{\lambda,k}^T \,\rose\SS =  \sfrac{1}{24}, \\
\sum_{\lambda=1}^s \Delta \rosc\SS[2]_\lambda\,(\rosc\SS_{\lambda-1} +  \rosc\SS_\lambda)\, \bar{\mathbf{r}}_{\lambda}^T\,(\boldsymbol{p}_\lambda + \boldsymbol{q}_\lambda) =  \sfrac{1}{12}, \\
\sum_{\lambda=1}^s \Delta \rosc\SS_\lambda\, \bar{\mathbf{r}}_{\lambda}^T\,\rosbeta\SS\,\rose\SS =  \sfrac{1}{24}, \\
%
\rosb\S* \,\sum_{\lambda=1}^s   \Delta \rosc\SS[2]_\lambda\,\ (\boldsymbol{p}_\lambda + \boldsymbol{q}_\lambda)\,(c\S_{\lambda}+2 \rosc\SS_{\lambda-1}) + \dots \qquad \\
\qquad \dots + 3\,\rosb\S* \,\sum_{\lambda=1}^s (\rosc\SS[2]_{\lambda} - \rosc\SS[2]_{\lambda-1}) \sum_{j=\lambda+1}^{s} \Delta \rosc\SS_j\, (\boldsymbol{p}_j + \boldsymbol{q}_j)  =  \sfrac{1}{4}, \\
\rosb\S* \,\sum_{\lambda=1}^s \Delta \rosc\SS_\lambda\, (\boldsymbol{p}_\lambda + \boldsymbol{q}_\lambda)\,\bar{\mathbf{r}}_{\lambda}^T\,\rose\SS =  \sfrac{1}{24}, \\
\rosb\S* \,\rosbeta\SS\,\sum_{\lambda=1}^s (\boldsymbol{p}_\lambda + \boldsymbol{q}_\lambda)\,(\rosc\SS[2]_{\lambda} - \rosc\SS[2]_{\lambda-1}) =  \sfrac{1}{12}.
\end{eqnarray*}
Note that, due to the structure of $\mathbf{p}_\lambda$, the second condition is a relation on the coefficients of the slow base scheme, and therefore it constrains the base methods that can be used.

\begin{remark}
Order four coupling conditions for infinitesimal step MR-GARK-ROW schemes can be derived in an analogous manner from the general order conditions \eqref{eqn:GARK-ROW-order4-conditions}.
\end{remark}
\else
\begin{remark}
Order four coupling conditions for infinitesimal step MR-GARK-ROS/ROW schemes can be derived in an analogous manner from the general order conditions  \eqref{eqn:GARK-ROS-order4-conditions}, \eqref{eqn:GARK-ROW-order4-conditions}.
\end{remark}
\fi

\section{Infinitesimal step SPC methods}
\label{sec:MRI-SPC-methods}
%
Consider a step-predictor-corrector method \eqref{eqn:MROS/MROW-SPC-butcher} with the slow base method  $(\rosb\S,\rosalpha\SS,\rosgamma\SS)$, and perform the compound step \eqref{eqn:MROS/MROW-SPC}.  Then proceed with the fast integration, but instead of applying the discrete formula \eqref{eqn:MROS/MROW-KF}, solve the following modified fast ODE:
\begin{eqnarray}
\label{eqn:MROS/MROW-SPC-inf}
\mathbf{v}' &=&\fun\F\Bigl( \mathbf{v} + \boldsymbol{\mu}^T({\scriptstyle \frac{\theta}{H}})\kron{\nvar} \mathbf{k}\S \Bigr)
 + \Lb\F \,\Bigl( \boldsymbol{\nu}^T({\scriptstyle \frac{\theta}{H}}) \kron{\nvar} \mathbf{k}\S \Bigr), \\
 \nonumber
&& \quad \mathbf{v}(0) = \y_{n-1}, \quad 0 \le \theta \le H.
\end{eqnarray}
The next step solution, computed using \eqref{eqn:MROS/MROW-Sol}, is:
\begin{equation}
\label{eqn:MROS/MROW-SPC-inf-solution}
\y_n = \mathbf{v}(H) +  \rosb\S*  \kron{\nvar} \mathbf{k}\S.
\end{equation}
\begin{remark}[Error estimation]
Assume that the base method has an embedded scheme $\hat{\rosb}\S$ for error estimation. Using $\hat{\rosb}\S$ in \eqref{eqn:MROS/MROW-SPC-inf-solution} gives an error estimate in the slow component, which can be used for macro-step error control. The fast components are solved with infinite accuracy, and the same fast integration \eqref{eqn:MROS/MROW-SPC-inf} is used for both the main and the embedded solutions.
\end{remark}

Solving \eqref{eqn:MROS/MROW-SPC-inf} with an arbitrarily accurate fast GARK-ROS/ROW method \eqref{eqn:GARK-ROS/ROW} with coefficients $(\rosb\F,\rosalpha\FF,\rosgamma\FF)$ gives the discrete solution:
\begin{equation}
\label{eqn:MROS/MROW-SPC-inf-fast}
\begin{split}
\mathbf{k}\F &= H\,\fun\F\left( 
\one_s \otimes \y_{n-1}  + \rosalpha\FF \kron{\nvar}  \mathbf{k}\F + \boldsymbol{\mu}^T(\rosc\FF)\kron{\nvar} \mathbf{k}\S \right)    \\
&\quad + \bigl( \Id_{s\F \times s\F} \otimes H\,\Lb\F \bigr)\,\left(\rosgamma\FF \kron{\nvar}  \mathbf{k}\F + \boldsymbol{\nu}^T(\rosc\FF) \kron{\nvar} \mathbf{k}\S \right),  \\ 
\mathbf{v}_{n} &= \y_{n-1} + \rosb\F*  \kron{\nvar}  \mathbf{k}\F.
\end{split}
\end{equation}

The Butcher tableau \eqref{eqn:mrGARK-ROW-butcher} of the coupled scheme \eqref{eqn:MROS/MROW-SPC} and \eqref{eqn:MROS/MROW-SPC-inf-fast} reads:
\begin{subequations}
\label{eqn:MROS/MROW-SPC-inf-butcher}
\begin{eqnarray}
\renewcommand{\arraystretch}{1.5}
\begin{array}{c|c|c}
\c\F & \boldsymbol{\upalpha}\FF & \boldsymbol{\upalpha}\FS \\ \hline
\c\S & \boldsymbol{\upalpha}\SF &\boldsymbol{\upalpha}\SS \\ \Xhline{1.5pt}
& \b\F* & \b\S*
\end{array} 
%
~\raisebox{-9pt}{=}~ 
\raisebox{9pt}{$
\begin{array}{c|cc|c}  
\rosc\SS & \rosalpha\SS & 0  & \rosalpha\SS \\
\rosc\FF & 0 & \rosalpha\FF    & \boldsymbol{\mu}^T(\rosc\FF) \\
\hline 
\rosc\SS & \rosalpha\SS & 0  & \rosalpha\SS   \\  
\Xhline{1.5pt}
 & 0 &  \rosb\F* &  \rosb\S* 
\end{array}
$},
\\
\begin{array}{c|c|c}  
\g\F & \boldsymbol{\upgamma}\FF & \boldsymbol{\upgamma}\FS \\ \hline
\g\S & \boldsymbol{\upgamma}\SF &\boldsymbol{\upgamma}\SS 
\end{array} 
~=~
\raisebox{5.5pt}{$
\begin{array}{c|cc|c}  
\rosg\SS & \rosgamma\SS & 0   & \rosgamma\SS \\
\rosg\FF &  0   & \rosgamma\FF & \boldsymbol{\nu}^T(\rosc\FF) \\
\hline 
\rosg\SS &\rosgamma\SS &  0 &  \rosgamma\SS 
 \end{array}
 $}.
\end{eqnarray}
\end{subequations}
The internal consistency conditions \eqref{eqn:ROW-internal-consistency} are:
\begin{equation}
\label{eqn:MROS/MROW-SPC-inf-internal-consistency}
\boldsymbol{\mu}^T(\rosc\FF)\,\one\F = \rosc \FF, \quad
\boldsymbol{\nu}^T(\rosc\FF)\,\one\F = \rosg \FF.
\end{equation}
Without loss of generality we choose $\rosgamma\FF = \zero$, i.e., we solve the modified fast ODE with an arbitrarily accurate Runge-Kutta scheme. The corresponding continuous coupling coefficients are chosen accordingly as $\boldsymbol{\nu}(t) = \zero \in \Re^{s\S \times 1}$.

Define the $\boldsymbol{\mu}(t) \in \Re^{s\S \times 1}$ coupling coefficients as polynomials in time:
\begin{equation*}
\begin{split}
&\boldsymbol{\mu}(t) \coloneqq \sum_{k \ge 0} \boldsymbol{\mu}_k \,t^k \quad \Rightarrow \quad \\
&\boldsymbol{\mu}^T(\rosc\FF) = \sum_{k \ge 0} \rosc \FF[\times k]\,\boldsymbol{\mu}_k^T, \qquad
\rosb\F* \,\boldsymbol{\mu}^T(\rosc\FF) = \sum_{k \ge 0} \frac{1}{k+1}\,\boldsymbol{\mu}_k^T \eqqcolon \overline{\boldsymbol{\mu}}^T. 
\end{split}
\end{equation*}
The internal consistency conditions \eqref{eqn:MROS/MROW-SPC-inf-internal-consistency} are satisfied with:
\[
 \boldsymbol{\mu}_1^T\,\one\S = 1; \quad \boldsymbol{\mu}_k^T\,\one\S = 0, ~ k\ne 1; \quad \boldsymbol{\nu}(t) = \zero_{s\S \times 1}. 
\]

\paragraph{SPC-MRI-GARK-ROS methods}
The order three ROS condition \eqref{eqn:GARK-ROS-order3-conditions} reads:
\begin{eqnarray*}
\overline{\boldsymbol{\mu}}^T\, \rose\SS = \sfrac{1}{6}.
\end{eqnarray*}
The order four ROS conditions \eqref{eqn:GARK-ROS-order4-conditions} are:
\ifnum \value{book}=1
\begin{eqnarray*}
\rosb\F* \, \Big((\boldsymbol{\mu}^T(\rosc\FF)\,\rose\SS) \times \rosc\FF \Big) = 
\sum_{k \ge 0} \frac{\boldsymbol{\mu}_k^T\,\rose\SS}{k+2} = \sfrac{1}{8}, \\
\rosb\F* \, \boldsymbol{\mu}^T(\rosc\FF)\,\rosc^{\times 2}  = 
\overline{\boldsymbol{\mu}}^T\,\rosc^{\times 2}  = \sfrac{1}{12}, \\ 
%
\rosb\F* \,\rosalpha\FF\,\boldsymbol{\mu}^T(\rosc\FF)\,\rose\SS = 
\sum_{k \ge 0} \left( \rosb\F* \,\rosalpha\FF\, \rosc\FF[\times k] \right)\,\boldsymbol{\mu}_k^T\,\rose\SS = \cdots \\
\qquad \cdots = 
\sum_{k \ge 0} \frac{\boldsymbol{\mu}_k^T\,\rose\SS}{(k+1)(k+2)} = \sfrac{1}{24}, \\
\overline{\boldsymbol{\mu}}^T\,\rosbeta\,\rose\SS = \sfrac{1}{24}.
%
\end{eqnarray*} 
\else
\begin{align*}
\sum_{k \ge 0} \frac{\boldsymbol{\mu}_k^T\,\rose\SS}{k+2} &= \sfrac{1}{8}, &
\overline{\boldsymbol{\mu}}^T\,\rosc\SS[\times 2] &= \sfrac{1}{12}, \\ 
\sum_{k \ge 0} \frac{\boldsymbol{\mu}_k^T\,\rose\SS}{(k+1)(k+2)} &= \sfrac{1}{24},&
\overline{\boldsymbol{\mu}}^T\,\rosbeta\SS\,\rose\SS &= \sfrac{1}{24}.
\end{align*} 
\fi
An order four linear-in-time coupling can be constructed as follows:
\begin{equation}
\label{eqn:ROS4-SPC-MRI-coupling}
\begin{split}
& \boldsymbol{\mu}(t) = \boldsymbol{\mu}_0 + \boldsymbol{\mu}_1\,t, \qquad \textnormal{where}\quad \boldsymbol{\mu}_0, \boldsymbol{\mu}_1\in \Re^{s\S \times 1}\quad \textnormal{satisfy:} \\
& \boldsymbol{\mu}_0^T\,\one\S = 0, \qquad \boldsymbol{\mu}_1^T\,\one\S = 1, \qquad
\boldsymbol{\mu}_0^T\,\rose\SS = -\sfrac{1}{12}, \qquad \boldsymbol{\mu}_1^T\,\rose\SS = \sfrac{1}{2}, \\
& \Big( \boldsymbol{\mu}_0^T + \sfrac{1}{2}\,\boldsymbol{\mu}_1^T \Big)\,\rosc^{\times 2}  = \sfrac{1}{12}, \qquad
\Big( \boldsymbol{\mu}_0^T + \sfrac{1}{2}\,\boldsymbol{\mu}_1^T \Big)\,\rosbeta\SS\,\rose\SS = \sfrac{1}{24}.
\end{split}
\end{equation}

\begin{example}[Multirate Rodas]
We build an SPC MRI version of Hairer and Wanner's Rodas method \cite[Chapter VI.4]{Hairer_book_II} by constructing a coupling of
the form \eqref{eqn:ROS4-SPC-MRI-coupling}. The Rodas method has six stages, and the coupling  \eqref{eqn:ROS4-SPC-MRI-coupling} is defined by twelve coefficients (the entries of the six-dimensional vectors $\boldsymbol{\mu}_0$ and $\boldsymbol{\mu}_1$). There are six ROS order four coupling conditions \eqref{eqn:ROS4-SPC-MRI-coupling}, and we also impose the order 3 ROW conditions \eqref{eqn:ROW3-SPC-MRI-order}. The coefficients depend on four free parameters:
\begin{equation}
\begin{split}
\boldsymbol{\mu}_0 & = 
\scalebox{0.7}{$
\begin{bmatrix}
\theta_1 \\
\theta_2 \\
- 1.923968128204745\,\theta_1 + 2.446324727549974\textnormal{e-01}\,\theta_2  + 4.509689603795104\textnormal{e-02} \\
1.405229246707428\,\theta_1 - 2.181782847233643\,\theta_2 + 1.289372580090594\textnormal{e-01} \\
- 4.812611185026834\textnormal{e-01}\,\theta_1 + 9.371503744786454\textnormal{e-01}\,\theta_2 - 1.740341540470105\textnormal{e-01} - \theta_3 \\
\theta_3
\end{bmatrix},
$}
\\
\boldsymbol{\mu}_1 &= 
\scalebox{0.7}{$
\begin{bmatrix}
- 2\,\theta_1   + 4.061438468864431\textnormal{e-01} \\
- 2\,\theta_2 + 5.932358823451654\textnormal{e-01} \\
3.847936256409489\,\theta_1 - 4.892649455099948\textnormal{e-01}\,\theta_2 - 3.657016798231872\textnormal{e-01} \\
- 2.810458493414855\,\theta_1 + 4.363565694467286\,\theta_2  + 5.003688760525202\textnormal{e-02} \\
9.625222370053667\textnormal{e-01}\,\theta_1 - 1.874300748957291\,\theta_2 + 3.162850629863266\textnormal{e-01} - \theta_4 \\
\theta_4
\end{bmatrix}.
$}
\end{split}
\end{equation}
The free parameters can be used to improve stability. By setting $\theta_3 = \theta_4 = 0$ the last slow Rodas stage is not used for coupling.
\end{example}

\paragraph{SPC-MRI-GARK-ROW methods}
\ifnum \value{book}=1
ROW order three conditions \eqref{eqn:GARK-ROW-order3-conditions}:
\begin{eqnarray*}
\rosb\F* \, \boldsymbol{\mu}^T(\rosc\FF)\, \rosc = \sfrac{1}{6}, \quad
\rosb\F* \, \boldsymbol{\mu}^T(\rosc\FF)\,\rosg = 0, \\
\rosb\F* \, \boldsymbol{\nu}^T(\rosc\FF)\, \rosc = 0, \quad
\rosb\F* \, \boldsymbol{\nu}^T(\rosc\FF)\,\rosg = 0.
\end{eqnarray*}
\fi
The order three ROW conditions \eqref{eqn:GARK-ROW-order3-conditions} are:
\begin{equation}
\label{eqn:ROW3-SPC-MRI-order}
\overline{\boldsymbol{\mu}}^T\, \rosc\SS = \sfrac{1}{6}, \qquad
\overline{\boldsymbol{\mu}}^T\,\rosg\SS = 0,
\end{equation}
and a third order coupling can be constructed as follows:
\begin{equation}
\label{eqn:ROW3-SPC-MRI-coupling}
\boldsymbol{\mu}(t) = \boldsymbol{\mu}_1\,t \qquad  \textnormal{with} \qquad
\boldsymbol{\mu}_1^T\,\one\S = 1, \quad 
\boldsymbol{\mu}_1^T\, \rosc\SS = \sfrac{1}{3}, \quad
\boldsymbol{\mu}_1^T\,\rosg\SS = 0.
\end{equation}
\begin{example}
Consider the third order, stiffly accurate ROW method ROS34PW2 of Rang and Angermann \cite[Table 4.3]{Rang_2005_W-methods}.
\ifnum \value{book}=1
The method coefficients are: 
\[
\scalebox{0.7}{$
\begin{array}{ll}
\rosgamma_{i,i} = 4.358665215084597\textnormal{e-01} & \\
\rosalpha_{2,1} =  8.7173304301691801\textnormal{e-01} & \rosgamma_{2,1} =  -8.7173304301691801\textnormal{e-01} \\
\rosalpha_{3,1} =  8.4457060015369423\textnormal{e-01} & \rosgamma_{3,1} =  -9.0338057013044082\textnormal{e-01} \\
\rosalpha_{3,2} =  -1.1299064236484185\textnormal{e-01} & \rosgamma_{3,2} =  5.4180672388095326\textnormal{e-02} \\
\rosalpha_{4,1} =  0 & \rosgamma_{4,1} =  2.4212380706095346\textnormal{e-01} \\
\rosalpha_{4,2} =  0 & \rosgamma_{4,2} =  -1.2232505839045147 \\
\rosalpha_{4,3} =  1 & \rosgamma_{4,3} =  5.4526025533510214\textnormal{e-01} \\
\rosb_{1} =  2.4212380706095346\textnormal{e-01} & \boldsymbol{\hat{b}}_{1} =  3.7810903145819369\textnormal{e-01} \\
\rosb_{2} =  -1.2232505839045147 & \boldsymbol{\hat{b}}_{2} =  -9.6042292212423178\textnormal{e-02} \\
\rosb_{3} =  1.5452602553351020 & \boldsymbol{\hat{b}}_{3} =  0.5 \\
\rosb_{4} =  4.3586652150845900\textnormal{e-01} & \boldsymbol{\hat{b}}_{4} =  2.1793326075422950\textnormal{e-01}
\end{array}
$}
\]
\fi
The SPC MRI coupling coefficients \eqref{eqn:ROW3-SPC-MRI-coupling} are defined in terms of one free parameter 
$\theta$ that can be used to improve stability:
\[
\scalebox{0.7}{$
\boldsymbol{\mu}_1 = 
\begin{bmatrix}
\theta \\
- 4.307016638790922 + 8.289196835086212\,\theta \\
4.541816529634874 - 7.686572272599903\,\theta \\
0.7652001091560487 - 1.602624562486310\,\theta
\end{bmatrix}.
$}
\]
\end{example}

\begin{remark}
Order four coupling conditions for infinitesimal SPC multirate GARK-ROW schemes can be derived in an analogous manner from the general order conditions \eqref{eqn:GARK-ROW-order4-conditions}.
\end{remark}

\section{Discussion}
\label{sec:conclusions}
This paper proposes a general framework for linearly-implicit multirate time integration.
Multirate GARK-ROS and GARK-ROW schemes, which make use of the exact or approximative Jacobian, respectively, are developed and analyzed.
Order conditions up to order four are derived, with a focus on internally consistent schemes.  We discuss several slow-fast coupling structures that lead to  efficient computational processes. Such couplings include compound-first step schemes, step-predictor-corrector methods, and multirate infinitesimal step approaches. Coefficient sets for new specific methods are given to illustrate these coupling strategies. 

The new MR(I)-GARK-ROS/ROW framework includes all existing multirate Rosenbrock(-W) methods as particular cases, and opens the possibility to develop new, high order, highly stable linearly implicit multirate schemes for a myriad of applications.  The development and optimization of practical MR(I)-GARK-ROS/ROW methods, their efficient implementation  \cite{Sandu_2014_FATODE}, and extensive numerical testing in real applications will be presented in a forthcoming publication.


\bibliographystyle{siam}

\end{document}